\documentclass[twocolumn,amssymb,aps,nofootinbib]{revtex4}
\usepackage{times,mathptmx,amsthm,amsmath}
\usepackage{pstricks,pst-node}
\makeatletter

\renewcommand{\p@subsection}{}

\renewcommand{\p@subsubsection}{}

\makeatother

\newtheorem{theorem}{Theorem}

\newtheorem{proposition}[theorem]{Proposition}
\newtheorem{corollary}[theorem]{Corollary}
\newtheorem{conjecture}[theorem]{Conjecture}
\newtheorem{example}[theorem]{Example}
\newtheorem{question}[theorem]{Question}
\theoremstyle{remark}
\newtheorem*{remark}{Remark}

\newcommand{\C}{\mathbb{C}}
\newcommand{\Z}{\mathbb{Z}}
\newcommand{\R}{\mathbb{R}}
\newcommand{\F}{\mathbb{F}}
\newcommand{\Sm}{\mathrm{Sm}}
\newcommand{\SL}{\mathrm{SL}}
\newcommand{\GL}{\mathrm{GL}}
\renewcommand{\tensor}{\otimes}
\DeclareMathOperator{\Pf}{Pf}
\DeclareMathOperator{\im}{im}
\DeclareMathOperator{\coker}{coker}
\newcommand{\braket}[1]{{\langle #1 \rangle}}

\newcommand{\vx}{\vec{x}}
\newcommand{\vq}{\vec{q}}
\newcommand{\tZ}{\tilde{Z}}
\newcommand{\tA}{\tilde{A}}

\newcommand{\eatline}{\vspace{-\baselineskip}}
\newcommand{\ie}{\textit{i.e.}}

\newcommand{\eq}[2]{\begin{equation}\label{#1}#2\end{equation}}

\newenvironment{fullfigure}[2]
    {\begin{figure}[htb]\begin{center}\def\ffa{#1}\def\ffb{#2}}
    {\vspace{\baselineskip}\caption{\ffb}\label{\ffa}\end{center}\end{figure}}

\newcommand{\fig}[1]{Figure~\ref{#1}}
\newcommand{\thm}[1]{Theorem~\ref{#1}}
\renewcommand{\sec}[1]{Section~\ref{#1}}

\psset{linewidth=.5pt,dash=4pt 4pt,unit=.25in,arrowsize=2pt 3}
\SpecialCoor

\newcommand{\mto}{\lput{:U}{\pspicture(0,0)(0,0)
\psline[arrows=->](2.3pt,0)(2.4pt,0)\endpspicture}}

\begin{document}
\title{Kasteleyn cokernels}
\author{Greg Kuperberg}
\thanks{Supported by NSF grants DMS \#9704125 and DMS \#0072342,
and by a Sloan Foundation Research Fellowship}
\affiliation{UC Davis}
\email{greg@math.ucdavis.edu}
\begin{abstract}
We consider Kasteleyn and Kasteleyn-Percus matrices, which arise in
enumerating matchings of planar graphs, up to matrix operations on their rows
and columns.  If such a matrix is defined over a principal ideal domain, this
is equivalent to considering its Smith normal form or its cokernel.  Many
variations of the enumeration methods result in equivalent matrices.  In
particular, Gessel-Viennot matrices are equivalent to Kasteleyn-Percus
matrices.

We apply these ideas to plane partitions and related planar of tilings. We
list a number of conjectures, supported by experiments in Maple, about the
forms of matrices associated to enumerations of plane partitions and other
lozenge tilings of planar regions and their symmetry classes.  We focus on the
case where the enumerations are round or $q$-round, and we conjecture that
cokernels remain round or $q$-round for related ``impossible enumerations'' in
which there are no tilings. Our conjectures provide a new view of the topic of
enumerating symmetry classes of plane partitions and their generalizations. In
particular we conjecture that a $q$-specialization of a Jacobi-Trudi matrix
has a Smith normal form.  If so it could be an interesting structure
associated to the corresponding irreducible representation of $\SL(n,\C)$.
Finally we find, with proof, the normal form of the matrix that appears in the
enumeration of domino tilings of an Aztec diamond.

\end{abstract}
\maketitle

\section{Introduction}
\label{s:intro}

The permanent-determinant and Hafnian-Pfaffian methods of Kasteleyn and Percus
give determinant and Pfaffian expressions for the number of perfect matchings
of a planar graph \cite{Kasteleyn:crystal,Percus:dimer}.  Although the methods
originated in mathematical physics, they have enjoyed new attention in
enumerative combinatorics in the past ten years
\cite{Jockusch:perfect,Kuperberg:perdet,Kuperberg:fun,Kenyon:local,Yang:thesis},
in particular for enumerating lozenge and domino tilings of various regions in
the plane.   These successes suggest looking at further properties of the
matrices that the methods produce beyond just their determinants or Pfaffians.

In this article we investigate the cokernel, or equivalently the Smith normal
form, of a Kasteleyn or Kasteleyn-Percus matrix $M$ arising from a planar graph
$G$.  One theme of our general results in Sections~\ref{s:gv} and \ref{s:invar} is that
the cokernel is a canonical object that can be defined in several different
ways.  More generally for weighted enumerations we consider $M$ up to the
equivalence relation of general row and column operations. If $G$ has at least
one matching, then the set of matchings is equinumerous with $\coker M$.  (In
\sec{s:bijective}, we conjecture an interpretation of this fact in the spirit
of a bijection.) The cokernel of $M$ is also interesting even when the graph
$G$ has no matchings, a situation which we call an \emph{impossible
enumeration}.  Propp proposed another invariant of $M$ that generalizes to
impossible enumerations and that was studied by Saldanha
\cite{Propp:progress,Saldanha:singular}, namely the spectrum of $M^*M$.

The idea of computing cokernels as a refinement of enumeration also arose in
the context of Kirchoff's determinant formula for the number of spanning trees
of a connected graph.  In this context the cokernels are called \emph{tree
groups} and they were proposed independently by Biggs, Lorenzini, and Merris
\cite{Biggs:algebraic,Lorenzini:attached,Merris:unimodular}. Indeed, Kenyon,
Propp, and Wilson \cite{KPW:trees}, generalizing an idea due to Fisher
\cite{Fisher:planar}, found a bijection between spanning trees of a certain
type of planar graph $G$ and the perfect matchings of another planar graph
$G'$.  We conjecture that the tree group of $G$ is isomorphic to the
Kasteleyn-Percus cokernel of $G'$.

In \sec{s:pp} we study cokernels for the special case of enumeration of plane
partitions in a box, as well as related lozenge tilings. We previously asked
what is the cokernel of a Carlitz matrix, which is equivalent to the
Kasteleyn-Percus matrix for plane partitions in a box with no symmetry imposed
\cite{Propp:progress}. We give two conjectures that together imply an
answer.  Finally in \sec{s:domino} we derive, with proof, the cokernel for
the enumeration of domino tilings of an Aztec diamond.

\acknowledgments

The author would like to thank Torsten Ekedahl, Christian Krattenthaler, and
Martin Loebl for helpful discussions.  The author would especially like to thank
Jim Propp for his diligent interest in this work.

\section{Preliminaries}
\label{s:prelim}

\subsection{Graph conventions}
\label{s:graph}

In general by a \emph{planar} graph we mean a graph embedded in
the sphere $S^2$.  We mark one point of $S^2$ outside of the graph
as the infinite point; the face containing it is the infinite face.
Our graphs may have both self-loops and multiple edges,
although self-loops cannot participate in matchings.

\subsection{Matrix algebra}
\label{s:matrix}

Let $R$ be a commutative ring with unit.  We consider matrices $M$ over
$R$, not necessarily square, up to three kinds of equivalence:  general row
operations,
$$M \mapsto AM$$
with $A$ invertible; general column operations,
$$M \mapsto MA$$
with $A$ invertible; and stabilization,
$$M \mapsto \left(\begin{array}{c|c}1 & 0 \\ \hline 0 & M \end{array}\right)$$
and its inverse.  Any matrix $M'$ which is equivalent to $M$
under these operations is a \emph{stably equivalent form} of $M$.

A matrix $A$ over $R$ is \emph{alternating} if it is antisymmetric and has null
diagonal.  (Antisymmetric implies alternating unless 2 is a zero divisor in
$R$.)  We consider alternating matrices up to two kinds of equivalence:
general symmetric operations,
$$A \mapsto BAB^T$$
with $B$ invertible; and stabilization,
$$A \mapsto \left(\begin{array}{c|c}
\begin{matrix} 0 & 1 \\ -1 & 0 \end{matrix} & 0 \\ \hline 0 & M \end{array}\right)$$
and its inverse.  A matrix $A'$ which is equivalent to $A$ is also called a
\emph{stably equivalent form} of $A$.

As a special case of these notions, \emph{elementary row operation} on a matrix
$M$ consists of either multiplying some row $i$ by a unit in $R$, or adding some
multiple of some row $i$ to row $j \ne i$.  \emph{Elementary column operations}
are defined likewise. We define a \emph{pivot} on a matrix $M$ at the $(i,j)$
position as subtracting $M_{k,j}/M_{i,j}$ times row $i$ from row $k$ for all $k
\ne i$, then subtracting $M_{i,k}/M_{i,j}$ times column $j$ from column $k$ for
all $j \ne k$.  This operation is possible when $M_{i,j}$ divides every entry
in the same row and column.  In matrix notation, if $M_{1,1} = 1$, the pivot at
$(1,1)$ looks like this:
$$M = \left(\begin{array}{c|c}1 & Y^T \\ \hline X & M' \end{array}\right)
    \mapsto \left(\begin{array}{c|c}1 & 0 \\ \hline 0 & M'-XY^T \end{array}
    \right).$$
A \emph{deleted pivot} consists of a pivot at $(i,j)$ followed by
deleting row $i$ and column $j$ from the matrix.  The deleted
pivot at $(1,1)$ on our example $M$ looks like this:
$$M = \left(\begin{array}{c|c}1 & Y^T \\ \hline X & M' \end{array}\right)
    \mapsto M'-XY^T.$$

If $A$ is an alternating matrix, we define an \emph{elementary symmetric
operation} as an elementary row operation followed by the same operation
in transpose on columns.  We can similarly define a \emph{symmetric
pivot} and a \emph{symmetric deleted pivot}.  All of these operations
are special cases of general symmetric matrix operations, and therefore
preserve the alternating property.

If $R$ is a principal ideal domain (PID), then an $n \times k$ matrix $M$ is
equivalent to one called its \emph{Smith normal form} and denoted
$\Sm(M)$.  We define $\Sm(M)$ and prove its existence in \sec{s:smith}.
Note that if $M'$ is a stabilization of $M$, then $\Sm(M')$
is a stabilization of $\Sm(M)$.

If $R$ is arbitrary, then we can interpret $M$ as a homomorphism
from $R^k$ to $R^n$.  In this interpretation $M$ has a kernel $\ker M$,
an image $\im M$, and a cokernel
$$\coker M = R^n/\im M.$$
If $R$ is a PID, the cokernel carries the same information as the Smith normal
form.  Over a general ring $R$, only very special matrices admit a Smith normal
form.  Determining equivalence of those that do not is much more complicated
than for those that do.  In particular inequivalent matrices may have the same
cokernel.  However, over any ring $R$ the cokernel is invariant under stable
equivalence and it does determine the determinant $\det M$ up to a unit
factor.  A special motivation for considering cokernels appears when $R = \Z$
and $M$ is square.  In this case the absolute determinant (\ie, absolute value
of the determinant) is the number of elements in the cokernel,
$$|\det M| = |\coker M|,$$
when the cokernel is finite, while
$$\det M = 0$$
if the cokernel is infinite.

An alternating matrix $A$ over a PID is also equivalent to its antisymmetric
Smith normal form $\Sm_a(A)$, which we also discuss in \sec{s:smith}.  Again
$\coker A$ determines $\Sm_a(A)$.

\begin{remark} If $M$ is a matrix over the polynomial ring $\F[x]$ over an
algebraically closed field $\F$, which is a PID, then the factor exhaustion
method for computing $\det M$ \cite{Krattenthaler:advanced} actually computes
the Smith normal form (or cokernel) of $M$.  Thus the Smith normal form plays a
hidden role in a computational method which is widely used in enumerative
combinatorics.

The basic version of the factor exhaustion method computes the rank of the
reduction
$$M \tensor \F[x]/(x-r)$$
for all $r \in \F$.  These ranks determine $\det M$ up to a constant factor if
the Smith normal form of $M$ is square free.  It is tempting to conclude that
the factor exhaustion method ``fails'' if the Smith normal form is not square
free. But sometimes one can compute the cokernel of
$$M \tensor \F[x]/(x-r)^k,$$
for all $r$ and $k$.  This information determines $\coker M$, as well as $\det
M$ up to a constant factor, regardless of its structure.  Thus the factor
exhaustion method always succeeds in principle.
\end{remark}

\section{Counting matchings}
\label{s:counting}

Most of this section is a review of Reference~\citealp{Kuperberg:fun}.

\subsection{Kasteleyn and Percus}
\label{s:perdet}

Let $G$ be a connected finite graph.  If we orient the edges of $G$, then we
define the \emph{alternating adjacency matrix} $A$ of $G$ by letting $A_{i,j}$
be the number of edges from vertex $i$ to vertex $j$ minus the number of edges
from vertex $j$ to vertex $i$.  If $G$ is simple, then the Pfaffian $\Pf A$ has
one non-zero term for every perfect matching of $G$, but in general the terms
may not have the same sign.

\begin{theorem}[Kasteleyn] If $G$ is a simple, planar graph, then it admits an
orientation such that all terms in $\Pf A$ have the same sign, where  $A$ is
the alternating adjacency matrix of $G$ \cite{Kasteleyn:crystal}.
\label{th:kasteleyn}
\end{theorem}

In general an orientation of $G$ such that all terms in $\Pf A$ have the same
sign is called a \emph{Pfaffian orientation} of $G$.  If an orientation of $G$
is Pfaffian, then the absolute Pfaffian $|\Pf A|$ is the number of perfect
matchings of $G$.  Kasteleyn's rule for a Pfaffian orientation is that an odd
number of edges of each (finite) face of $G$ should point clockwise.  We call
such an orientation \emph{Kasteleyn flat} and the resulting matrix $A$ a
\emph{Kasteleyn matrix} for the graph $G$.  Likewise an orientation may be
Kasteleyn flat at a particular face if it satisfies Kasteleyn's rule at that
face. Every planar graph has a Kasteleyn-flat orientation, although it is only
flat at the infinite face of $G$ if $G$ has an even number of vertices.
Forming a Kasteleyn matrix to count matchings of a planar graph is also called
the \emph{Hafnian-Pfaffian method} \cite{Kuperberg:perdet}.

Percus \cite{Percus:dimer} found a simplification of the Hafnian-Pfaffian
method when $G$ is bipartite.  Suppose that $G$ is a bipartite graph with the
vertices colored black and white, and suppose that each edge has a sign $+$ or
$-$, interpreted as the weight $1$ or $-1$. Then we define the \emph{bipartite
adjacency matrix} $M$ of $G$ by letting $M_{i,j}$ be the total weight of all
edges from the black vertex $i$ to the white vertex $j$.  If $G$ is simple,
then the determinant $\det M$ has a non-zero term for each perfect matching of
$G$, but in general with both signs.

\begin{theorem}[Percus] If $G$ is a simple, planar, bipartite graph, then it
admits a sign decoration such that all terms in $\det M$ have the same sign,
where $M$ is the bipartite adjacency matrix of  $G$.
\label{th:percus}
\end{theorem}

In the rule given by Percus, the edges of each face of $G$ should
have an odd number of $-$ signs if and only if the face has $4k$ sides.
We call such a sign decoration of $G$ \emph{Kasteleyn flat} and the
corresponding matrix $M$ a \emph{Kasteleyn-Percus matrix}. Every planar graph
has a Kasteleyn-flat signing, although it is only flat at the infinite
face if $G$ has an even number of vertices.  Forming a Kasteleyn-Percus matrix
$M$ is also called the \emph{permanent-determinant method}. A Kasteleyn matrix
$A$ for a bipartite graph $G$ can be viewed as two copies of a Kasteleyn-Percus
matrix $M$:
$$A = \left(\begin{array}{c|c} 0 & M \\ \hline -M & 0 \end{array}\right).$$

\begin{fullfigure}{f:triple}{Tripling an edge in a graph.}
\pspicture(-1,-.5)(6.5,.5)
\psline(-.3,-.3)(0,0)(-.3,.3)\psline(-.5,0)(1.5,0)
\psline(1.3,.3)(1,0)(1.3,-.3)
\qdisk(0,0){.1}
\qdisk(1,0){.1}
\psline{->}(2,0)(3,0)
\psline(3.7,-.3)(4,0)(3.7,.3)
\psline(3.5,0)(6.5,0)
\psline(6.3,.3)(6,0)(6.3,-.3)
\qdisk(4,0){.1}
\qdisk(4.67,0){.1}
\qdisk(5.33,0){.1}
\qdisk(6,0){.1}
\endpspicture
\end{fullfigure}

If the graph $G$ is not simple, then we may make it simple by tripling edges,
as shown in \fig{f:triple}. The set of matchings of the new graph $G'$ is
naturally bijective with the set of matchings of $G$.  A more economical
approach is to define a  Kasteleyn matrix $A$ or a Kasteleyn-Percus matrix $M$
for $G$ directly. In this case $A_{ij}$ is the number of edges from vertex $i$
to vertex $j$ minus the number from $j$ to $i$, while $M_{ij}$ is the number of
positive edges minus the number of negative edges connecting $i$ and $j$.

A variant of the Hafnian-Pfaffian method applies to a projectively planar graph
$G$ which is \emph{locally but not globally bipartite}. This means that $G$ is
embedded in the projective plane and that all faces have an even number of
sides, but that $G$ is not bipartite.  An equivalent condition is that all
contractible cycles in $G$ have even length and all non-contractible cycles
have odd length.

\begin{theorem} If a projectively planar graph $G$ is locally but not globally
bipartite, then it admits a Pfaffian orientation. \cite{Kuperberg:fun}
\label{th:proj}
\end{theorem}

The orientation constructed in the proof of Theorem~\ref{th:proj} is one with
the property that each face has an odd number of edges pointing in each
direction.   We call such an orientation Kasteleyn flat; it exists if $G$ has
an even number of vertices.  (If $G$ has an odd number of vertices, then every
orientation is trivially Pfaffian.)  We call the corresponding alternating
adjacency matrix $A$ the Kasteleyn matrix of $G$ as usual.

The constructions of this section, in particular Theorems~\ref{th:kasteleyn}
and \ref{th:percus}, generalize to weighted enumerations of the matchings of
$G$, where each edge of $G$ is assigned a weight and the weight of a matching
is the product of the weights of its edges.  We separately assign signs or
orientations to $G$ using the Kasteleyn rule (in the general case) or the
Percus rule (in the bipartite case).  The weighted alternating adjacency matrix
$A$ is called a Kasteleyn matrix of $G$.  If $G$ is bipartite, the weighted
bipartite adjacency matrix $M$, with the weights multiplied by the signs, is a
Kasteleyn-Percus matrix of $G$.  Then $\Pf A$ or $\det M$ is the total weight
of all matchings in $G$.

\subsection{Polygamy and reflections}
\label{s:polygamy}

We can use the Hafnian-Pfaffian method to count certain generalized matchings
among the vertices of a planar graph $G$ using an idea originally due to Fisher
\cite{Fisher:planar}. We arbitrarily divide the vertices of $G$ into three
types:  Monogamous vertices, odd-polygamous vertices, and even-polygamous
vertices.  An \emph{odd-polygamous vertex} is one that can be connected to any
odd number of other vertices in a matching, while an \emph{even-polygamous
vertex} can be connected to any even number of other vertices (including
none).

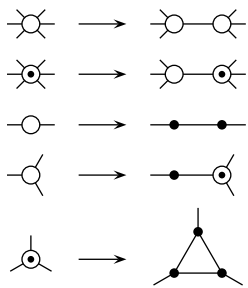
\begin{fullfigure}{f:polyres}{Resolving polygamy into monogamy.}
\pspicture(-1,-.5)(5,.5)
\psline(-.3,-.3)(0,0)(-.3,.3)\psline(-.5,0)(.5,0)\psline(.3,.3)(0,0)(.3,-.3)
\pscircle[fillstyle=solid](0,0){.2}
\psline{->}(1,0)(2,0)
\psline(2.7,-.3)(3,0)(2.7,.3)\psline(2.5,0)(4.5,0)\psline(4.3,.3)(4,0)(4.3,-.3)
\pscircle[fillstyle=solid](3,0){.2}
\pscircle[fillstyle=solid](4,0){.2}
\endpspicture \\
\pspicture(-1,-.5)(5,.5)
\psline(-.3,-.3)(0,0)(-.3,.3)\psline(-.5,0)(.5,0)\psline(.3,.3)(0,0)(.3,-.3)
\pscircle[fillstyle=solid](0,0){.2}\qdisk(0,0){.07}
\psline{->}(1,0)(2,0)
\psline(2.7,-.3)(3,0)(2.7,.3)\psline(2.5,0)(4.5,0)\psline(4.3,.3)(4,0)(4.3,-.3)
\pscircle[fillstyle=solid](3,0){.2}
\pscircle[fillstyle=solid](4,0){.2}\qdisk(4,0){.07}
\endpspicture \\
\pspicture(-1,-.5)(5,.5)
\psline(-.5,0)(.5,0)
\pscircle[fillstyle=solid](0,0){.2}
\psline{->}(1,0)(2,0)
\psline(2.5,0)(4.5,0)
\qdisk(3,0){.1}\qdisk(4,0){.1}
\endpspicture \\
\pspicture(-1,-.5)(5,.5)
\psline(-.5,0)(0,0)(.25,.433)\psline(0,0)(.25,-.433)
\pscircle[fillstyle=solid](0,0){.2}
\psline{->}(1,0)(2,0)
\psline(2.5,0)(4,0)(4.25,.433)\psline(4,0)(4.25,-.433)
\qdisk(3,0){.1}\pscircle[fillstyle=solid](4,0){.2}\qdisk(4,0){.07}
\endpspicture \\
\pspicture(-1,-.5)(5,1.5)
\psline(0,.789)(0,.289)(.433,.039)
\psline(0,.289)(-.433,.039)
\pscircle[fillstyle=solid](0,.289){.2}\qdisk(0,.289){.07}
\psline{->}(1,.289)(2,.289)
\pspolygon(3,0)(4,0)(3.5,.866)(3,0)
\psline(3,0)(2.567,-.25)\psline(4,0)(4.433,-.25)\psline(3.5,.866)(3.5,1.366)
\qdisk(3,0){.1}\qdisk(4,0){.1}\qdisk(3.5,.866){.1}
\endpspicture
\end{fullfigure}

If $G$ is a graph with polygamous vertices, we can find a new graph $G'$ such
that the ordinary perfect matchings of $G'$ are bijective with the generalized
matchings of $G$.  The graph $G'$ defined from $G$ using a series of local
moves that are shown in \fig{f:polyres}.  (In this figure and later, an open
circle is an even-polygamous vertex and a dotted circle is an odd-polygamous
vertex.)  We also describe the moves in words.  First, if a polygamous vertex
of $G$ has valence greater than 3, we can split it into two polygamous vertices
of lower valence with the same total parity. This leaves polygamous vertices of
both parities of valence 1, 2, and 3.  If a polygamous vertex $v$ is even and
has valence 1, we can delete it.  If it is odd and has valence 1 or 2, it is
the same as an ordinary monogamous vertex.  If it is even and has valence 2, we
can replace it with two monogamous vertices.  If it is even and has valence 3,
we can split it into an odd-polygamous divalent vertex and an odd-polygamous
trivalent vertex.  Finally, if it odd and has valence 3, we can replace it with
a triangle.  Each of these moves comes with an obvious bijection between the
matchings before and after.  Thus these moves establish the following:

\begin{proposition}[Fisher]
Given a graph $G$ with odd- and even-polygamous vertices, the polygamous
vertices can be replaced by monogamous subgraphs so that the matchings of the
new graph $G'$ are bijective with those of $G$.  If $G$ is planar,
then $G'$ can be planar.
\end{proposition}

We call the resulting graph $G'$ a \emph{monogamous resolution} of $G$.
If $G$ is planar, then $G'$ admits Kasteleyn matrices, and we call
any such matrix a Kasteleyn matrix of $G$ as well.

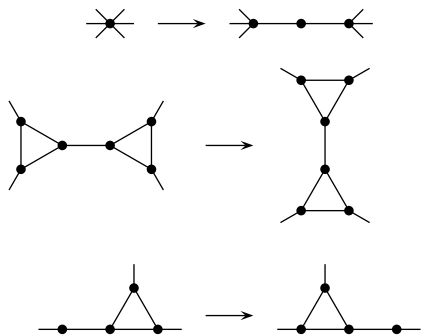
\begin{fullfigure}{f:resmoves}{Moves on monogamous resolutions
of a polygamous graph.}
\pspicture(-1,-.5)(5.5,.5)
\psline(-.3,-.3)(.3,.3)
\psline(-.3,.3)(.3,-.3)
\psline(-.5,0)(.5,0)
\qdisk(0,0){.1}
\psline{->}(1,0)(2,0)
\psline(2.7,-.3)(3,0)(2.7,.3)
\psline(2.5,0)(5.5,0)
\psline(5.3,.3)(5,0)(5.3,-.3)
\qdisk(3,0){.1}
\qdisk(4,0){.1}
\qdisk(5,0){.1}
\endpspicture \\
\pspicture(-1,-2)(6.5,2)
\psline(-.5,0)(-1.366,.5)(-1.366,-.5)(-.5,0)
    (.5,0)(1.366,.5)(1.366,-.5)(.5,0)
\qdisk(-.5,0){.1}\qdisk(.5,0){.1}
\psline(-1.366,.5)(-1.616,.933)\psline(1.366,.5)(1.616,.933)
\psline(-1.366,-.5)(-1.616,-.933)\psline(1.366,-.5)(1.616,-.933)
\qdisk(-1.366,.5){.1}\qdisk(1.366,.5){.1}
\qdisk(-1.366,-.5){.1}\qdisk(1.366,-.5){.1}
\psline{->}(2.5,0)(3.5,0)
\psline(5,-.5)(5.5,-1.366)(4.5,-1.366)(5,-.5)
    (5,.5)(5.5,1.366)(4.5,1.366)(5,.5)
\qdisk(5,-.5){.1}\qdisk(5,.5){.1}
\psline(5.5,-1.366)(5.933,-1.616)\psline(5.5,1.366)(5.933,1.616)
\psline(4.5,-1.366)(4.067,-1.616)\psline(4.5,1.366)(4.067,1.616)
\qdisk(5.5,-1.366){.1}\qdisk(5.5,1.366){.1}
\qdisk(4.5,-1.366){.1}\qdisk(4.5,1.366){.1}
\endpspicture \\
\pspicture(-2,-.5)(5.5,1.8)
\psline(-2,0)(1,0)\psline(.5,0)(0,.866)(-.5,0)
\psline(0,.866)(0,1.366)
\qdisk(-.5,0){.1}\qdisk(.5,0){.1}
\qdisk(0,.866){.1}\qdisk(-1.5,0){.1}
\psline{->}(1.5,.289)(2.5,.289)
\psline(6,0)(3,0)\psline(3.5,0)(4,.866)(4.5,0)
\psline(4,.866)(4,1.366)
\qdisk(4.5,0){.1}\qdisk(3.5,0){.1}
\qdisk(4,.866){.1}\qdisk(5.5,0){.1}
\endpspicture
\end{fullfigure}

The monogamous resolution of a polygamous graph is far from unique.   But we
can consider moves that connect different monogamous resolutions of a
polygamous graph.   The moves are as shown in \fig{f:resmoves}: Doubly
splitting a vertex, rotating a pair of triangles, and switching a triangle with
an edge.  Each of these moves comes with a bijection between the matchings of
the two graphs that it connects.

\begin{proposition} Any two monogamous resolutions of a graph $G$
are connected by the moves of vertex splitting and its inverse,
switching triangles, and switching a triangle with an edge.
The moves also connect any two planar resolutions of a planar graph $G$
through intermediate planar resolutions.
\label{p:resmoves}
\end{proposition}

The proof of Proposition~\ref{p:resmoves} is routine.

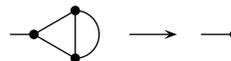
\begin{fullfigure}{f:selfconn}{Removing a self-connected triangle.}
\pspicture(-.5,-1)(4.5,1)
\psline(-.5,0)(0,0)(.866,-.5)(.866,.5)(0,0)
\qdisk(0,0){.1}
\qdisk(.866,-.5){.1}
\qdisk(.866,.5){.1}
\psarc(.866,0){.5}{-90}{90}
\psline{->}(2,0)(3,0)
\psline(3.5,0)(4.2,0)
\qdisk(4.2,0){.1}
\endpspicture\eatline
\end{fullfigure}

Another interesting move is removing a self-connected triangle,
as shown in \fig{f:selfconn}.
This move induces a 2-to-1 map on the set of matchings before and after.

Polygamous matchings have two common applications.  If a graph $G$ is entirely
polygamous, then we can denote the presence or absence of each edge by an
element of $\Z/2$.  Each vertex then imposes a linear constraint on the
variables, so the number of matchings is therefore either 0 or $2^n$ for some
$n$.  The corresponding weighted enumerations are related to the Ising model
\cite{Kuperberg:fun,Waerden:reichweite,Fisher:planar}.  Another way to see that
the number of matchings is a power of two is to use the moves in
Figures~\ref{f:resmoves} and \ref{f:selfconn} to reduce a monogamous
resolution of $G$ to a tree, which has at most one matching.

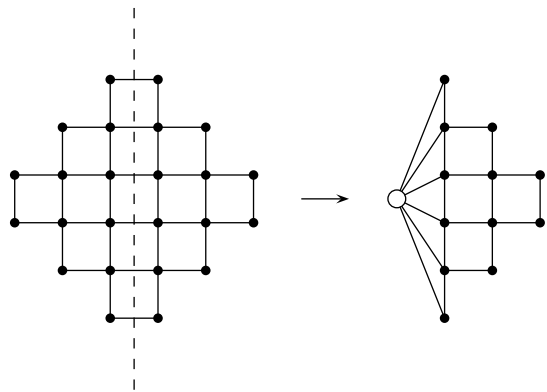
\begin{fullfigure}{f:polytied}{Using polygamy to count reflection-invariant
    matchings.}
\pspicture(-3,-4)(9,4)
\pspolygon(-2.5,.5)(-2.5,-.5)(2.5,-.5)(2.5,.5)
\pspolygon(-1.5,1.5)(-1.5,-1.5)(1.5,-1.5)(1.5,1.5)
\pspolygon(.5,-2.5)(-.5,-2.5)(-.5,2.5)(.5,2.5)
\multips( -.5, 2.5)(1,0){2}{\qdisk(0,0){.1}}
\multips(-1.5, 1.5)(1,0){4}{\qdisk(0,0){.1}}
\multips(-2.5,  .5)(1,0){6}{\qdisk(0,0){.1}}
\multips(-2.5, -.5)(1,0){6}{\qdisk(0,0){.1}}
\multips(-1.5,-1.5)(1,0){4}{\qdisk(0,0){.1}}
\multips( -.5,-2.5)(1,0){2}{\qdisk(0,0){.1}}
\psline[linestyle=dashed](0,-4)(0,4)
\psline{->}(3.5,0)(4.5,0)
\pspolygon(5.5,0)(6.5,2.5)(6.5,-2.5)
\pspolygon(5.5,0)(6.5,1.5)(7.5,1.5)(7.5,-1.5)(6.5,-1.5)
\pspolygon(5.5,0)(6.5,.5)(8.5,.5)(8.5,-.5)(6.5,-.5)
\multips(6.5,2.5)(0,-1){6}{\qdisk(0,0){.1}}
\multips(7.5,1.5)(0,-1){4}{\qdisk(0,0){.1}}
\multips(8.5, .5)(0,-1){2}{\qdisk(0,0){.1}}
\pscircle[fillstyle=solid](5.5,0){.2}
\endpspicture
\end{fullfigure}

Another application is counting matchings invariant under reflections
\cite{Kuperberg:perdet,Kuperberg:fun}.  Suppose that a planar graph $G$ has a
reflection symmetry $\sigma$, and suppose that the line of reflection bisects
some of the edges of $G$.  Then the $\sigma$-invariant matchings of $G$ are
bijective with a modified quotient graph $G//\sigma$ in which the bisected
edges are tied to a polygamous vertex, as in \fig{f:polytied}.
The parity of the polygamous vertex should be set so that the total parity
(odd-polygamous plus monogamous vertices) is even. The same construction works
if we divide $G$ by any group acting on the sphere that includes reflections,
since all of the reflective boundary can be reached by a single polygamous
vertex.

\subsection{Gessel-Viennot}
\label{s:gv}

The Gessel-Viennot method \cite{GV:binomial,GV:partitions} yields another
determinant expression for a certain sum over the sets of disjoint paths in an
acyclic, directed graph $G$.  (Theorem~\ref{th:gv} below, which is the basic
result of the method, was independently found by Lindstr\"om
\cite{Lindstrom:vector}.  Gessel and Viennot were the first to use it for
unweighted enumeration.) The graph $G$ need not be planar.  We label some of
the vertices of $G$ as left endpoints and some as right endpoints, and we
separately order the left endpoints and the right endpoints. Let $P$ be the set
of collections of vertex-disjoint paths in $G$ connecting the left endpoints to
the right endpoints. If $P$ is non-empty then there are the same number of left
and right endpoints on left and right; if there are $n$ of each we call the
elements of $P$ \emph{disjoint $n$-paths}. The Gessel-Viennot matrix $V$ is
defined by setting $V_{i,j}$ to the number of paths in $G$ from left endpoint
$i$ to right endpoint $j$.

\begin{theorem}[Lindstr\"om, Gessel-Viennot]
Let $G$ be a directed, acyclic, weighted graph with $n$ ordered left endpoints
and $n$ ordered right endpoints.  Let $P$ be the set of disjoint $n$-paths in
$G$ connecting left to right.  If $V$ is the Gessel-Viennot matrix of $G$, then
\eq{e:gv}{\det V = \sum_{\ell \in P} w(\ell) (-1)^\ell.}
Here $(-1)^\ell$ is the sign of the bijection from the left to the right
endpoints induced by the paths in the collection $\ell$, and $w(\ell)$ is the
product of the weights of the edges of $G$ that appear in $\ell$.
\label{th:gv}
\end{theorem}
\begin{proof} We outline a non-traditional proof that will be useful later.  We
first suppose that the left endpoints are the sources in $G$ (the vertices with
in-degree 0) and the right endpoints are the sinks (the vertices with
out-degree 0). We argue by induction on the number of \emph{transit
vertices}, meaning vertices that are neither sources nor sinks.

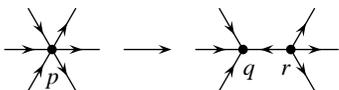
\begin{fullfigure}{f:gvsplit}{Splitting a transit vertex.}
\pspicture(-1,-1)(6,1)
\pcline(-.5,.866)(0,0)\mto\pcline(-1,0)(0,0)\mto\pcline(-.5,-.866)(0,0)\mto
\pcline(0,0)(1,0)\mto\pcline(0,0)(.5,.866)\mto\pcline(0,0)(.5,-.866)\mto
\qdisk(0,0){.1}
\rput[t](0,-.5){$p$}
\psline{->}(1.5,0)(2.5,0)
\pcline(3.5,.866)(4,0)\mto\pcline(3,0)(4,0)\mto\pcline(3.5,-.866)(4,0)\mto
\pcline(5,0)(6,0)\mto\pcline(5,0)(5.5,.866)\mto\pcline(5,0)(5.5,-.866)\mto
\pcline(5,0)(4,0)\mto\qdisk(4,0){.1}\qdisk(5,0){.1}
\rput[tl](4,-.3){$q$}
\rput[tr](5,-.3){$r$}
\endpspicture\eatline
\end{fullfigure}

If $G$ has no transit vertices, every path in $G$ has length one. Consequently
the $n$-paths in $G$ are perfect matchings, and equation~\eqref{e:gv}
is equivalent to the definition of the determinant.  Suppose then that
$p$ is a transit vertex in $G$.  We form a new graph $G'$ by splitting $p$
into two vertices $q$ and $r$, with $q$ a sink and $r$ a source,
as shown in \fig{f:gvsplit}.
We number $q$ and $r$ as the $n+1$st (last) source and sink in $G'$. We give
the new edge between $q$ and $r$ a weight of $-1$. There is a natural bijection
between disjoint $n$-paths $\ell$ in $G$ and disjoint $n+1$-paths $\ell'$ in
$G'$:  Every path in $\ell$ which avoids $p$ is included in $\ell'$.  If some
path in $\ell$ meets $p$, we break it into two paths ending at $q$ and starting
again at $r$.  If $\ell$ is disjoint from $p$, we include the edge from $r$ to
$q$ in $\ell'$.  In order to argue that the right side of
equation~\eqref{e:gv} are the same for $G$ and $G'$,
we check that
$$(-1)^\ell w(\ell) = (-1)^{\ell'} w(\ell').$$
If $\ell$ avoids $p$, the two sides are immediately the same. If $\ell$ meets
$p$, then $(-1)^\ell$ and $(-1)^{\ell'}$ have opposite sign and so do $w(\ell)$
and $w(\ell')$.  The left side of equation~\eqref{e:gv} is also the same:  If
$V$ and $V'$ are the Gessel-Viennot matrices of $G$ and $G'$, $V$ is obtained
from $V'$ by a deleted pivot at $(n+1,n+1)$.

Now suppose that the left and right endpoints do not coincide with the sources
and sinks.  If $G$ has a left endpoint $q$ which is not a source, then there is
an edge $e$ from a vertex $p$ to the vertex $q$.  Let $G'$ be $G$ with $e$
removed and let $V'$ be its Gessel-Viennot matrix.   Since the edge $e$ is not
in any $n$-path in $G$, the graph $G'$ has the same $n$-paths with the same
weights.  If $p$ is the $i$th left endpoint, we can obtain $V'$ from $V$ by
subtracting $w(i,j)$ times row $i$ from $j$, where $w(i,j)$ is the total weight
in $G$ of all paths from $i$ to $j$.  These row operations do not change the
determinant.  The same argument applies if $G$ has a right endpoint which is
not a sink.

Finally if $G$ has a right endpoint source or a left endpoint sink,
then it has no $n$-paths and some row or column of $V$ is 0.
If $G$ has a source or a sink which is not an endpoint,
we can delete it without changing the Gessel-Viennot matrix $V$ or
the set of $n$-paths.
\end{proof}

We call the graph $G'$ constructed in our proof of \thm{th:gv} the
\emph{transit-free resolution of $G$}.  The transit-free resolution is a
connection between the Gessel-Viennot method and the permanent-determinant
method:

\begin{corollary} Let $G$ be a connected, planar, directed, acyclic graph with
$n$ left and right endpoints on the outside face.  Suppose that the left
endpoints are segregated from the right endpoints on this face. If the left
endpoints are the sources and the right endpoints are the sinks, then the
Gessel-Viennot matrix $V$ of $G$ is obtained from a Kasteleyn-Percus matrix $M$
of the transit-free resolution $G'$ of $G$ by deleted pivots.  If not every
left endpoint is a source or not every right endpoint is a sink, $V$ is
obtained by deleted pivots and other matrix operations.
\label{c:gv}
\end{corollary}

Note that by construction the matchings of the transit-free resolution $G'$ of
$G$ are bijective with the $n$-paths in $G$.  The planarity of $G$ together
with the position of its endpoints imply that all $n$-paths induce the same
bijection and therefore have the same sign.

\begin{proof} The Gessel-Viennot matrix $V'$ of the graph $G'$ is obtained from
$V$ by the stated operations.  Thus it suffices to show that $G'$ is planar and
that $V'$ is also a  Kasteleyn-Percus matrix of $G'$.

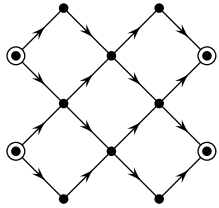
\begin{fullfigure}{f:lr}{Left-to-right orientation implied by segregation
    of sources and sinks (circled).}
\pspicture(-2,-2)(2,2)
\pcline(-2, 1)(-1, 2)\mto\pcline(-2, 1)(-1, 0)\mto
\pcline(-2,-1)(-1, 0)\mto\pcline(-2,-1)(-1,-2)\mto
\pcline(-1, 2)( 0, 1)\mto\pcline(-1, 0)( 0, 1)\mto
\pcline(-1, 0)( 0,-1)\mto\pcline(-1,-2)( 0,-1)\mto
\pcline( 0, 1)( 1, 2)\mto\pcline( 0, 1)( 1, 0)\mto
\pcline( 0,-1)( 1, 0)\mto\pcline( 0,-1)( 1,-2)\mto
\pcline( 1, 2)( 2, 1)\mto\pcline( 1, 0)( 2, 1)\mto
\pcline( 1, 0)( 2,-1)\mto\pcline( 1,-2)( 2,-1)\mto
\pscircle[fillstyle=solid](-2,1){.2}
\pscircle[fillstyle=solid](-2,-1){.2}
\pscircle[fillstyle=solid]( 2,1){.2}
\pscircle[fillstyle=solid]( 2,-1){.2}
\qdisk(-2, 1){.1}\qdisk(-2,-1){.1}
\qdisk(-1, 2){.1}\qdisk(-1, 0){.1}\qdisk(-1,-2){.1}
\qdisk( 0, 1){.1}\qdisk( 0,-1){.1}
\qdisk( 1, 2){.1}\qdisk( 1, 0){.1}\qdisk( 1,-2){.1}
\qdisk( 2, 1){.1}\qdisk( 2,-1){.1}
\endpspicture
\end{fullfigure}

We first establish that the orientation of $G$ is qualitatively like that of
the example in \fig{f:lr}: the orientations all point from left to right.
More precisely, the edges incident to each transit vertex are segregated, in
the sense that all incoming edges are adjacent and all outgoing edges are
adjacent.  The edges of each internal face are also segregated, in the sense
that the clockwise edges are adjacent and the counterclockwise edges are
adjacent.  To prove that $G$ has this structure, we reallocate the Euler
characteristic of the sphere, $2$, expressed as a sum over elements of $G$.
In this sum, each vertex and face has Euler characteristic $1$ and each edge
has Euler characteristic $-1$.  If a pair of edges shares both a vertex $v$
and a face $f$, we deduct $\frac12$ from the Euler characteristic of $f$ if
the edges both point to or both point from $v$, and otherwise we deduct
$\frac12$ from $v$.  Since each edge participates in $4$ such pairs, these
deductions absorb the total Euler characteristic of all edges.

The reallocated characteristic of a vertex is $1$ if it is a source or sink,
$0$ if it is a segregated transit vertex, and negative otherwise.  The
reallocated characteristic of a face is at most $2-2n$ if it is the outside
face (since orientations must switch between clockwise and counterclockwise
both at the sources and sinks and between them), $0$ if it is a segregated
internal face, and negative otherwise.  (No face has positive reallocated
characteristic since $G$ is acyclic.)  Thus the only way that the total can be
$2$ is if all internal faces and all transit vertices are segregated.

That the transit vertices of $G$ are segregated implies that $G'$ is planar.
That each internal face $f$ is segregated implies that if $f$ has $k$ sides,
the corresponding face $f'$ of $G'$ has $2k-2$ sides.  Moreover the $k-2$ new
edges of $f'$ have weight $-1$ in the proof of \thm{th:gv}, which agrees with
the Kasteleyn-Percus rule.  Thus $V'$ is a Kasteleyn-Percus matrix of $G'$, as
desired.
\end{proof}

Finally we have not discussed a Pfaffian version of the Gessel-Viennot method
defined by Stembridge \cite{Stembridge:paths}. We believe that this method can
be generalized further, and that  it admits an analogue of
Corollary~\ref{c:gv}.

\section{Matchings and Smith normal form}
\label{s:snf}

\subsection{Equivalences of Kasteleyn and Kasteleyn-Percus matrices}
\label{s:invar}

If $M$ is a Kasteleyn-Percus matrix of a bipartite, planar graph $G$, then we
can consider its cokernel, which by \sec{s:matrix} is equinumerous with the
number of matchings of $M$ if it has any matchings. Furthermore, if $\coker M$
is infinite or if $M$ isn't square, we can think of $\coker M$ as a way to
``count'' matchings in a graph that has none. We call such a computation
an \emph{impossible enumeration}. Both observations are reasons to study
$\coker M$ as part of enumerative combinatorics.

If $G$ is weighted by elements of some ring $R$, then we can consider $M$ up to
stable equivalence, whether or not it has a Smith normal form.

Suppose that $M$ and $M'$ are two Kasteleyn-Percus matrices for the same
planar graph $G$.  Then the signs on $G$ given by $M'$ and $M$ differ by a
1-cocycle $c$ with coefficients in the group $\{+,-\}$ \cite{Kuperberg:fun}.
Since the sphere has no first homology, $c = \delta d$, where $d$ is a
0-cochain.  More explicitly $d$ is a function from the vertices of $G$ to
$\{+,-\}$.   We can use $d$ to form two diagonal matrices $A$ and $B$ with
diagonal entries $\pm 1$ and such that $M' = AMB$. Evidently $A$ and $B$ are
invertible over $\Z$, so $M'$ and $M$ have the same cokernel.  In conclusion:

\begin{proposition}  If $G$ be a weighted bipartite planar graph, then all of
its Kasteleyn-Percus matrices $M$ are stably equivalent forms.  In particular
$\coker M$ is an invariant of $G$.
\label{p:inv}
\end{proposition}

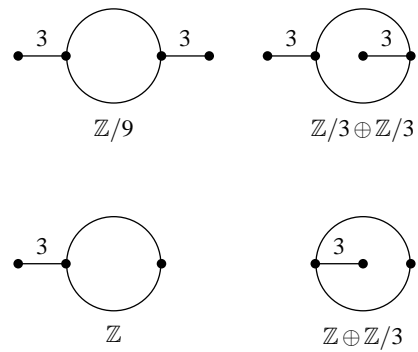
\begin{fullfigure}{f:embed}{Embedding-dependent Kasteleyn-Percus cokernels.}
\begin{tabular}{cc}
\pspicture(-2.5,-2)(2.5,1)
\psline(-2,0)(-1,0)\rput[b](-1.5,.2){3}
\pscircle(0,0){1}
\psline(1,0)(2,0)\rput[b](1.5,.2){3}
\qdisk(-2,0){.1}\qdisk(-1,0){.1}\qdisk(1,0){.1}\qdisk(2,0){.1}
\rput[t](0,-1.3){$\Z/9$}
\endpspicture &
\pspicture(-2.5,-2)(2.5,1)
\psline(-2,0)(-1,0)\rput[b](-1.5,.2){3}
\pscircle(0,0){1}
\psline(0,0)(1,0)\rput[b](.5,.2){3}
\qdisk(-2,0){.1}\qdisk(-1,0){.1}\qdisk(0,0){.1}\qdisk(1,0){.1}
\rput[t](0,-1.3){$\Z/3 \oplus \Z/3$}
\endpspicture \\[.75cm]
\pspicture(-2.5,-2)(2.5,1)
\psline(-2,0)(-1,0)\rput[b](-1.5,.2){3}
\pscircle(0,0){1}
\qdisk(-2,0){.1}\qdisk(-1,0){.1}\qdisk(1,0){.1}
\rput[t](0,-1.3){$\Z$}
\endpspicture &
\pspicture(-2.5,-2)(2.5,1)
\psline(-1,0)(0,0)\rput[b](-.5,.2){3}
\pscircle(0,0){1}
\qdisk(-1,0){.1}\qdisk(0,0){.1}\qdisk(1,0){.1}
\rput[t](0,-1.3){$\Z \oplus \Z/3$}
\endpspicture\eatline
\end{tabular}
\end{fullfigure}

\begin{example}  The Smith normal form or cokernel of $M$ can depend on the
embedding of $G$ in the plane.  The top two graphs in \fig{f:embed} have
cokernels $\Z/9$ and $\Z/3 \oplus \Z/3$ in the two embeddings shown.  If $G$
has an odd number of vertices, then it cannot be Kasteleyn flat on its outside
face. In this case changing which face is on the outside can change the
cokernel as well.  The bottom two graphs in \fig{f:embed} are an example.
\end{example}


Our analysis generalizes to the non-bipartite case.  If $G$ is
a planar graph with a Kasteleyn matrix $A$, then we can consider
$A$ up to equivalence.

Again all Kasteleyn matrices we choose for the planar graph $G$ are equivalent,
because any two Kasteleyn-flat orientations of $G$ differ by the coboundary of
a 0-cochain on $G$ with values in $\{+,-\}$.  The matrices are consequently
equivalent under the transformation
$$A \mapsto B^TAB$$
for some diagonal matrix $B$ whose non-zero entries are $\pm 1$.  We can also
pass from the usual clockwise-odd Kasteleyn rule to the counterclockwise-odd
rule by negating $A$.  We have no reason to believe that $A$ and $-A$ are
equivalent over a general ground ring $R$, but they do have the same
cokernel and are therefore equivalent if $R$ is a PID.

If $G$ is projectively planar and locally but not
globally bipartite, the argument is slightly different while the conclusion is
the same.  In this case the cohomology group
$$H^1(\R P^2,\Z/2) \cong \Z/2$$
is non-trivial.  Suppose that we have two Kasteleyn-flat orientations of $G$
whose matrices are $A$ and $A'$.  Their discrepancy is a 1-cocycle $c$ which
could represent either the trivial or the non-trivial class in $H^1(\R
P^2,\Z/2)$.  If $c$ is trivial, then
$$A' = B^TAB$$
for some diagonal $B$.  I.e., $A$ and $A'$ are equivalent.  If $c$ is
non-trivial, then
$$A' = -B^TAB,$$
i.e., $A'$ is equivalent to $-A$.

\begin{fullfigure}{f:doublevert}{Preserving Kasteleyn flatness.}
\pspicture(0,-.5)(6.5,.5)
\psline(.7,-.3)(1,0)(.7,.3)\psline(.5,0)(1.5,0)
\psline(1.3,.3)(1,0)(1.3,-.3)
\qdisk(1,0){.1}
\psline{->}(2,0)(3,0)
\psline(3.7,-.3)(4,0)(3.7,.3)
\psline(3.5,0)(4,0)\psline(6,0)(6.5,0)
\pcline(4,0)(5,0)\mto\pcline(6,0)(5,0)\mto
\psline(6.3,.3)(6,0)(6.3,-.3)
\qdisk(4,0){.1}\qdisk(5,0){.1}\qdisk(6,0){.1}
\endpspicture\eatline
\end{fullfigure}

Kasteleyn and Kasteleyn-Percus matrices remain equivalent under more operations than
just the choices of signs or orientations. In particular they remain equivalent
under the moves in Section~\ref{s:polygamy}.  In each move we make a graph $G'$
from the graph $G$, and we need to choose related Kasteleyn-flat orientations
of both graphs.  For example consider a double vertex splitting. If $G$ is
Kasteleyn-flat, and if we orient the two new edges in the splitting in opposite
directions as in \fig{f:doublevert},  then $G'$ is also Kasteleyn flat.  If the
three vertices are numbered $1$, $2$, and $3$, then the matrix $A'$ of $G'$ has
a submatrix of the form
$$\begin{pmatrix} 0 & 1 & 0 \\ -1 & 0 & -1 \\ 0 & 1 & 0 \end{pmatrix}.$$
If we perform a deleted pivot at $(1,2)$ and $(2,1)$,
we reduce $A'$ to the matrix $A$ of $G$.

\subsection{Is it bijective?}
\label{s:bijective}

Whenever two sets are known to have the same size, a traditional question in
combinatorics is whether or not there is a bijection between them.  In this
section we conjecture a relationship between cokernels and matching sets
which is similar to a bijection.

If $M$ is a Kasteleyn-Percus matrix for an unweighted bipartite, planar graph
$G$ with at least one perfect matching, then two such sets to consider are
$\coker M$ and $P$, the set of perfect matchings of
$G$. To this pair we must add a
third set, $\coker M^T$, since the choice between $M$ and $M^T$ is
arbitrary. As explained in \sec{s:smith}, $\coker M$ and $\coker M^T$ are
isomorphic, but there is no canonical isomorphism.  This is evidence against a
natural bijection between $\coker M$ and $\coker M^T$, and therefore a natural
bijection between either of them and $P$.  On the other hand, the special
planar structure of $M$ might yield such bijections.

It may be better to consider quantum bijections or linearized bijections
instead of traditional ones.  If $A$ and $B$ are two finite sets,
a \emph{quantum bijection} is a unitary isomorphism
$$\C[A] \cong \C[B]$$
between the formal linear spans of $A$ and $B$. A quantum bijection can be
implemented by a quantum computer algorithm just as a traditional bijection can
be implemented on a standard computer \cite{Steane:review}.  A linear bijection is a linear
isomorphism
$$\F[A] \cong \F[B],$$
not necessarily unitary, for some field $\F$.  A linear bijection does not have
the empirical computational interpretation that a traditional bijection or a
quantum bijection does, but as a means of proving that $A$ and $B$ are
equinumerous, it can be considered constructive.

If $M$ is a non-singular $n \times n$ matrix over $\Z$, then there is a natural
quantum bijection between $\coker M$ and $\coker M^T$, namely the discrete
Fourier transform.  (It is also a special case of Pontryagin duality
\cite{Rudin:fourier}).  We express it by defining a unitary matrix $U$ whose
rows are indexed by $x \in \coker M$ and whose columns are indexed by $y \in
\coker M^T$. Given such $x$ and $y$, we let $X$ and $Y$ be lifts in $\Z^n$. We
then define
$$U_{x,y} = \frac{\exp(2\pi i Y^T M^{-1} X)}{\sqrt{|\det M|}}.$$
It is easy to check that $Y^T M^{-1} X$ changes by an integer if we change the
lift $X$ of $x$, because two such lifts differ by an element in $\im M$.

Given $M$, $G$, and $P$ as above, it might be possible to factor the unitary
map $U$ into maps to and from $\C[P]$.  However we may need to further relax
the notion of a bijection.   Sometimes when $G$ is a finite group equinumerous
with a finite set $S$, there is no natural bijection between them, but instead
there is a natural, freely transitive group action of $G$ on $S$.  Having
introduced quantum bijections, we can try to make $\C[P]$ a free unitary module
over the group algebras  $\C[\coker M]$ and $\C[\coker M^T]$.  We can even ask
that the two group actions be compatible with $U$ by requiring the commutation
relations
$$\alpha_y \alpha_x = \exp(2\pi i Y^T M^{-1} X)\alpha_x \alpha_y,$$
where $x \in \coker M$ and $y \in \coker M^T$ and $\alpha_x$ and $\alpha_y$ are
their hypothetical actions on $\C[P]$.  A standard theorem in representation
theory says that for any $M$ the algebra
$$D = \C[\coker M] \tensor \C[\coker M^T]$$
twisted by this commutation relation is
isomorphic to a matrix algebra.  This means that $D$ has only one irreducible
representation, and we conjecture that $\C[P]$ has the structure
of this representation.

If $A$ is a non-singular alternating matrix, then we can define a similar
algebra $D$ as a deformation of the group algebra $\C[\coker A]$. We let $D$ be
the formal complex span of elements $\alpha_x$ with $x \in \coker A$, and
we arbitrarily order the elements in $\coker A$.  For $x \le y \in \coker A$,
we define
$$\alpha_x \alpha_y = \alpha_{x+y},$$
and for arbitrary $x$ and $y$, we impose the relation
$$\alpha_y \alpha_x = \exp(2\pi i Y^T A^{-1} X) \alpha_x \alpha_y.$$
The algebra $D$ is again isomorphic to a matrix algebra, since
\thm{th:asmith} allows us to put $A$ into
$$A = \left(\begin{array}{c|c} 0 & M \\ \hline -M & 0 \end{array}\right).$$
If we do so then algebra $D$ then has the form previously described.

\begin{conjecture} If $A$ is a Kasteleyn matrix of a planar graph $G$
with non-empty matching set $P$, then there is a natural action of
the algebra $D$ on $\C[P]$, possibly depending on the way that $G$
is embedded in the plane.
\label{c:action}
\end{conjecture}

\section{Lozenge and domino tilings}
\label{s:lozdom}

\subsection{Plane partitions and lozenge tilings}
\label{s:pp}

Plane partitions are an interesting source of enumerative planar matching
problems.  We consider the cokernels and Smith normal forms that arise in these
problems, relying on the material in References~\citealp{Kuperberg:perdet} and
\citealp{Kuperberg:fun}.

\begin{fullfigure}{f:pptil}{A plane partition and a tiling.}
\psset{xunit=.5cm,yunit=.577cm}
\pspicture(-7,0)(7,5)
\rput[b](-4.5,0){\pspicture(-2,0)(2,5)
\psset{fillstyle=solid}
\pspolygon[fillcolor=lightgray](-1,4.5)(-1,3.5)(0,3)(0,4)
\pspolygon[fillcolor=darkgray](1,4.5)(1,3.5)(0,3)(0,4)
\pspolygon[fillcolor=lightgray](-1,3.5)(-1,2.5)(0,2)(0,3)
\pspolygon[fillcolor=lightgray](0,3)(0,2)(1,1.5)(1,2.5)
\pspolygon[fillcolor=darkgray](2,3)(2,2)(1,1.5)(1,2.5)
\pspolygon[fillcolor=lightgray](0,2)(0,1)(1,.5)(1,1.5)
\pspolygon[fillcolor=darkgray](2,2)(2,1)(1,.5)(1,1.5)
\pspolygon[fillcolor=lightgray](-2,2)(-2,1)(-1,.5)(-1,1.5)
\pspolygon[fillcolor=darkgray](0,2)(0,1)(-1,.5)(-1,1.5)
\pspolygon[fillcolor=gray](-1,4.5)(0,4)(1,4.5)(0,5)
\pspolygon[fillcolor=gray](0,3)(1,2.5)(2,3)(1,3.5)
\pspolygon[fillcolor=gray](-2,2)(-1,1.5)(0,2)(-1,2.5)
\psset{fillstyle=none}
\psline[linestyle=dashed](-2,2)(-2,4)(-1,4.5)
\psline[linestyle=dashed](1,4.5)(2,4)(2,3)
\psline[linestyle=dashed](-1,.5)(0,0)(1,.5)(0,1)
\pcline[linestyle=none](-2,1)(-2,4)\Aput{$c$}
\pcline[linestyle=none](0,0)(-2,1)\Aput{$b$}
\pcline[linestyle=none](0,0)(2,1)\Bput{$a$}
\endpspicture}
\rput[b](4.5,0){\pspicture(-2,0)(2,5)
\pspolygon(-1,4.5)(-1,3.5)(0,3)(0,4)
\pspolygon(1,4.5)(1,3.5)(0,3)(0,4)
\pspolygon(-1,3.5)(-1,2.5)(0,2)(0,3)
\pspolygon(0,3)(0,2)(1,1.5)(1,2.5)
\pspolygon(2,3)(2,2)(1,1.5)(1,2.5)
\pspolygon(0,2)(0,1)(1,.5)(1,1.5)
\pspolygon(2,2)(2,1)(1,.5)(1,1.5)
\pspolygon(-2,2)(-2,1)(-1,.5)(-1,1.5)
\pspolygon(0,2)(0,1)(-1,.5)(-1,1.5)
\pspolygon(-1,4.5)(0,4)(1,4.5)(0,5)
\pspolygon(0,3)(1,2.5)(2,3)(1,3.5)
\pspolygon(-2,2)(-1,1.5)(0,2)(-1,2.5)
\psline(-2,2)(-2,4)(-1,4.5)
\psline(1,4.5)(2,4)(2,3)
\psline(-1,.5)(0,0)(1,.5)(0,1)
\psline(-2,3)(-1,3.5)
\pcline[linestyle=none](-2,1)(-2,4)\Aput{$c$}
\pcline[linestyle=none](0,0)(-2,1)\Aput{$b$}
\pcline[linestyle=none](0,0)(2,1)\Bput{$a$}
\endpspicture}
\psline{->}(-1,2.5)(1,2.5)
\endpspicture
\end{fullfigure}

\begin{fullfigure}{f:z223}{The graph $Z(2,2,3)$.}
\psset{xunit=.289cm,yunit=.5cm}
\pspicture(-6,1)(6,9)
\multips(-4,2)(0,2){4}{\qline(0,0)(2,0)}  
\multips(-1,1)(0,2){5}{\qline(0,0)(2,0)}
\multips( 2,2)(0,2){4}{\qline(0,0)(2,0)}
\multips(-4,2)(0,2){3}{\qline(0,0)(-1,1)} 
\multips(-1,1)(0,2){4}{\qline(0,0)(-1,1)}
\multips( 2,2)(0,2){4}{\qline(0,0)(-1,1)}
\multips( 5,3)(0,2){3}{\qline(0,0)(-1,1)}
\multips( 4,2)(0,2){3}{\qline(0,0)( 1,1)} 
\multips( 1,1)(0,2){4}{\qline(0,0)( 1,1)}
\multips(-2,2)(0,2){4}{\qline(0,0)( 1,1)}
\multips(-5,3)(0,2){3}{\qline(0,0)( 1,1)}
\multips(-5,3)(0,2){3}{\qdisk(0,0){.1}}\multips(-2,2)(0,2){4}{\qdisk(0,0){.1}}
\multips( 1,1)(0,2){5}{\qdisk(0,0){.1}}\multips( 4,2)(0,2){4}{\qdisk(0,0){.1}}
\multips( 5,3)(0,2){3}{\qdisk(0,0){.1}}\multips( 2,2)(0,2){4}{\qdisk(0,0){.1}}
\multips(-1,1)(0,2){5}{\qdisk(0,0){.1}}\multips(-4,2)(0,2){4}{\qdisk(0,0){.1}}
\endpspicture
\end{fullfigure}

A plane partition in an $a \times b \times c$ box is equivalent to a lozenge
tiling of an $(a,b,c)$-semiregular hexagon (\fig{f:pptil}), which in
turn is equivalent to a perfect matching of the dual hexagonal graph
$Z(a,b,c)$ (\fig{f:z223}).

If a group $G$ acts on the box and the plane partitions inside it, it also acts
on the graph $Z(a,b,c)$.  The $G$-invariant matchings correspond to matchings
of a modified quotient graph $Z_G(a,b,c)$.  To understand these
graphs we recall the three generators of $G$:
\begin{description}
\item[$\bullet$] $\rho$, cyclic symmetry for plane partitions or
rotation by 120 degrees for lozenge tilings, defined when $a=b=c$.
\item[$\bullet$] $\tau$, symmetry for plane partitions or diagonal
reflection for lozenge tilings, defined when $b=c$.
\item[$\bullet$] $\kappa$, complementation for plane partitions or
rotation by 180 degrees for lozenge tilings.
\end{description}
We describe $Z_G(a,b,c)$ case by case:
\begin{description}
\item[$\bullet$] If $G = \braket{\rho}$ or $\braket{\rho,\kappa}$, or if $a$,
$b$, and $c$ are all even and $G$ is $\braket{\kappa}$, then $Z_G(a,b,c)$
is the usual quotient graph $Z(a,b,c)/G$.
\item[$\bullet$] If $G = \braket{\kappa\tau}$ or $G =
\braket{\kappa\tau,\rho}$, then we delete the edges and vertices of $Z(a,b,c)$
along the lines of reflection and let $Z_G(a,b,c)$ be a connected component of
the remainder.
\item[$\bullet$] If $G = \braket{\kappa}$ and only one or two of $a$, $b$, and
$c$ is even,  then $Z(a,b,c)$ has a central edge $e$ invariant under $\kappa$.
If one dimension is even, we define $Z_\kappa(a,b,c)$ by deleting $e$ and its
vertices and then quotienting by $\kappa$. If two dimensions are even, we
define $Z_\kappa(a,b,c)$ by deleting $e$ but not its vertices, and then
quotienting by $\kappa$.
\item[$\bullet$] If $G = \braket{\tau}$, $\braket{\tau,\kappa}$, or
$\braket{\tau,\rho,\kappa}$, we start with $Z_H(a,b,c)$, where $H$ is,
respectively, $\braket{1}$, $\braket{\tau\kappa}$, or
$\braket{\rho,\tau\kappa}$. We cut $Z_H(a,b,c)$ by the line of reflection, and
tie the cut edges of one region to a polygamous vertex to define $Z_G(a,b,c)$.
\item[$\bullet$] If $G = \braket{\tau,\rho}$, we cut $Z(a,a,a)$ by three lines
of reflection and tie the cut edges of one region (lying along two of the
three lines) to a polygamous vertex to form $Z_G(a,a,a)$.
\end{description}

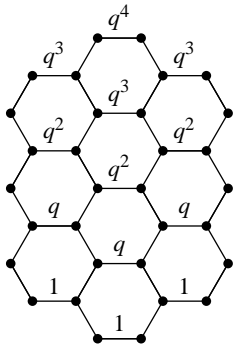
\begin{fullfigure}{f:qz223}{The graph $Z(2,2,3)$ weighted for $q$-enumeration.}
\psset{xunit=.289cm,yunit=.5cm}
\pspicture(-6,1)(6,9)
\multips(-4,2)(0,2){4}{\qline(0,0)(2,0)}  
\multips(-1,1)(0,2){5}{\qline(0,0)(2,0)}
\multips( 2,2)(0,2){4}{\qline(0,0)(2,0)}
\multips(-4,2)(0,2){3}{\qline(0,0)(-1,1)} 
\multips(-1,1)(0,2){4}{\qline(0,0)(-1,1)}
\multips( 2,2)(0,2){4}{\qline(0,0)(-1,1)}
\multips( 5,3)(0,2){3}{\qline(0,0)(-1,1)}
\multips( 4,2)(0,2){3}{\qline(0,0)( 1,1)} 
\multips( 1,1)(0,2){4}{\qline(0,0)( 1,1)}
\multips(-2,2)(0,2){4}{\qline(0,0)( 1,1)}
\multips(-5,3)(0,2){3}{\qline(0,0)( 1,1)}
\qline(-1,1)( 1,1) \qline(-5,3)(-4,2)
\qline(-2,2)(-1,3) \qline( 1,3)( 2,2)
\qline( 4,2)( 5,3) \qline(-4,4)(-2,4)
\qline( 1,5)( 2,4) \qline( 4,4)( 5,5)
\qline(-5,5)(-4,6) \qline(-2,6)(-1,5)
\qline( 2,6)( 4,6) \qline(-5,7)(-4,8)
\qline(-2,8)(-1,7) \qline( 1,7)( 2,8)
\qline( 4,8)( 5,7) \qline(-1,9)( 1,9)
\multips(-5,3)(0,2){3}{\qdisk(0,0){.1}}\multips(-2,2)(0,2){4}{\qdisk(0,0){.1}}
\multips( 1,1)(0,2){5}{\qdisk(0,0){.1}}\multips( 4,2)(0,2){4}{\qdisk(0,0){.1}}
\multips( 5,3)(0,2){3}{\qdisk(0,0){.1}}\multips( 2,2)(0,2){4}{\qdisk(0,0){.1}}
\multips(-1,1)(0,2){5}{\qdisk(0,0){.1}}\multips(-4,2)(0,2){4}{\qdisk(0,0){.1}}
\rput[b](-3,8.2){$q^3$} \rput[b](-3,6.2){$q^2$} \rput[b](-3,4.2){$q$}
\rput[b](-3,2.2){$1$}   \rput[b](0,9.2){$q^4$}  \rput[b](0,7.2){$q^3$}
\rput[b](0,5.2){$q^2$}  \rput[b](0,3.2){$q$}    \rput[b](0,1.2){$1$}
\rput[b](3,8.2){$q^3$}  \rput[b](3,6.2){$q^2$}  \rput[b](3,4.2){$q$}
\rput[b](3,2.2){$1$}
\endpspicture
\end{fullfigure}

For every symmetry group $G$, the number of $G$-invariant plane partitions
$N_G(a,b,c)$ is \emph{round}, meaning a product of small factors.  (This is not
a completely rigorous notion.  In this case $N_G(a,b,c)$ grows exponentially in
$ab+ac+bc$ while its prime factors grow linearly in $a+b+c$.) In addition if
$G$ is a subgroup of $\braket{\rho,\tau}$, the $G$-invariant partitions have
round $q$-enumerations, where the $q$-weight of a plane partition is sometimes
the number of cubes and sometimes the number of $G$-orbits of cubes.  (To say
that a polynomial $P(q)$ is round means not only that it is a product of small
factors, but also that the factors are cyclotomic.  Equivalently $P(q)$ is a
ratio of products of differences of monomials.)  All of these enumerations are
proven \cite{Stanley:symmetries,Kuperberg:perdet,%
Andrews:tsscpp,Stembridge:enumeration} except for the conjectured orbit
$q$-enumeration of totally symmetric plane partitions.  Both orbit and cube
$q$-enumerations can be realized by weighted enumerations of matchings in
$Z_G(a,b,c)$ \cite{Kuperberg:perdet,Kuperberg:fun}.  \fig{f:qz223} shows an
example where $G$ is trivial (so that there is no distinction between cube and
orbit $q$-enumeration in this case).   Also the weight of a matching in the
example only agrees with the $q$-weight of the corresponding plane partition up
to a (matching-independent) factor of $q$.  We will therefore consider
$q$-enumerations over the ring $\Z[q,q^{-1}]$ and absorb powers of $q$ in
normalization.

We let $Z_G(a,b,c;q)$ be $Z_G(a,b,c)$ with weights chosen for
cube $q$-enumeration, and we let $\tZ_G(a,b,c;q)$ be $Z_G(a,b,c)$
with weights chosen for orbit $q$-enumeration.

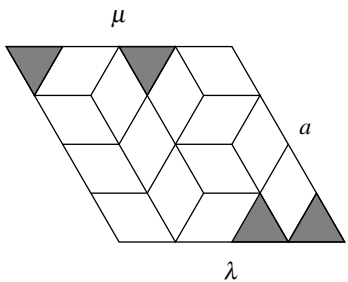
\begin{fullfigure}{f:lmregion}{A lozenge tiling of a region with $a=4$,
    $\lambda=(2,2)$, and $\mu=(1)$.}
\psset{xunit=.75cm,yunit=.650cm}
\pspicture(-2,0)(4,5)
\pspolygon[fillstyle=solid,fillcolor=gray](-2,4)(-1,4)(-1.5,3)
\pspolygon[fillstyle=solid,fillcolor=gray](0,4)(1,4)(.5,3)
\pspolygon[fillstyle=solid,fillcolor=gray](2,0)(2.5,1)(3,0)
\pspolygon[fillstyle=solid,fillcolor=gray](3,0)(3.5,1)(4,0)
\psline(-1,4)(-1.5,3)(-.5,3)(0,4)(1.5,1)
\psline(-1,2)(0,2)(1,4)(1.5,3)(2.5,3)(2,2)
\psline(-.5,1)(.5,1)(1.5,3)
\psline(-.5,3)(1,0)(1.5,1)(2.5,1)(3,2)
\psline(1,2)(2,2)(3,0)(3.5,1)
\pspolygon(-2,4)(2,4)(4,0)(0,0)
\rput[t](2,-.4){$\lambda$}
\rput[b](0,4.4){$\mu$}
\rput[lb](3.2,2.2){$a$}
\endpspicture
\end{fullfigure}

\begin{fullfigure}{f:lmweights}{Weights of $Z(\lambda/\mu;\vx)$.}
\psset{xunit=1cm,yunit=.577cm}
\pspicture(-2,.5)(4,6)
\psline(-1,4)(-.5,3.5)
\psline(-.5,2.5)(0,2)
\multips(-1,4)(1,0){4}{\psline(0,0)(0,1)\qdisk(0,0){.1}\qdisk(0,1){.1}}
\multips(-.5,2.5)(1,0){4}{\psline(0,0)(0,1)\qdisk(0,0){.1}\qdisk(0,1){.1}}
\multips(0,1)(1,0){4}{\psline(0,0)(0,1)\qdisk(0,0){.1}\qdisk(0,1){.1}}
\multips(-1,5)(2,0){2}{\psline(0,0)(.5,.5)(1,0)\qdisk(.5,.5){.1}}
\multips(-.5,3.5)(1,0){3}{\psline(0,0)(.5,.5)(1,0)}
\multips(0,2)(1,0){3}{\psline(0,0)(.5,.5)(1,0)}
\multips(0,1)(1,0){2}{\psline(0,0)(.5,-.5)(1,0)\qdisk(.5,-.5){.1}}
\rput[rb](-.8,5.3){$x_4$}\rput[rb](1.2,5.3){$x_4$}
\rput[rb](-.3,3.8){$x_3$}\rput[rb](.7,3.8){$x_3$}\rput[rb](1.7,3.8){$x_3$}
\rput[rb](.2,2.3){$x_2$}\rput[rb](1.2,2.3){$x_2$}\rput[rb](2.2,2.3){$x_2$}
\rput[rb](.7,.8){$x_1$}\rput[rb](1.7,.8){$x_1$}
\endpspicture
\end{fullfigure}

The problem of counting matchings in $Z(a,b,c)$ has two other interesting
generalizations.  The lozenge tilings of a parallelogram strip with notches on
both sides (see \fig{f:lmregion}) are naturally bijective with the
semi-standard skew tableaux of shape $\lambda/\mu$ and parts bounded by some
$a$ \cite[\S7.10]{Stanley:enumerative2}, where $\lambda$ and $\mu$ are two
partitions with $\lambda$ containing $\mu$.  If $\lambda$ has $b$ nonzero
parts, the parallelogram then has dimensions $a$ by $\lambda_1 + b$.  On the
bottom row it has notches at positions $\lambda_i + b+1-i$, counting from the
left.  On the top row it has notches at positions $\mu_i + b+1-i$, counting
from the left and extending $\mu$ by $0$ so that it also has $b$ parts.  We
leave the proof of the bijection between tilings and tableaux to the reader
since it is similar to existing arguments in the literature.  We let
$Z(\lambda/\mu,a)$ be the graph dual to the tiling of this region by triangles,
so that the lozenge tilings correspond to the matchings of $Z(\lambda/\mu,a)$.
In an important weighted enumeration we assign the weight $x_i$ to the
northeast-pointing edges in the $i$th row, as in \fig{f:lmweights}. Call the
weighted graph $Z(\lambda/\mu;\vx)$ The total weight of its matchings is the
skew Schur function $s_{\lambda/\mu}(\vx)$ (in finitely many variables).  If we
let
$$\vq_a = (1,q,\ldots,q^{a-1}),$$
then the specialization $s_{\lambda/\mu}(\vq_a)$ is the standard $q$-enumeration
of skew tableaux.  In particular, if we omit $\mu$ (by setting it to the empty
partition), then $s_\lambda(\vx)$ is the character of an irreducible
representation $V(\lambda)$ of $\GL(a,\C)$, while $s_\lambda(\vq_a)$ is round by
the $q$-Weyl dimension formula.

A second generalization is to count tilings of a semiregular hexagon with side
lengths $a,b+d,c,a+d,b,c+d$ and with a triangle of size $d$ removed
\cite{CEKZ:triangular}.  We do not know of a $q$-enumeration of these, although
the symmetry $\kappa\tau$ (corresponding to TCPPs) appears when $b=c$ and the
symmetry $\rho$ (corresponding to CSPPs) appears when $a=b=c$.  We let
$Z(a,b,c,d,d)$ be the dual graph and we analogously
define $Z_G(a,b,c,d,d)$ if $G$ is a subgroup of
$\braket{\kappa\tau,\rho}$.

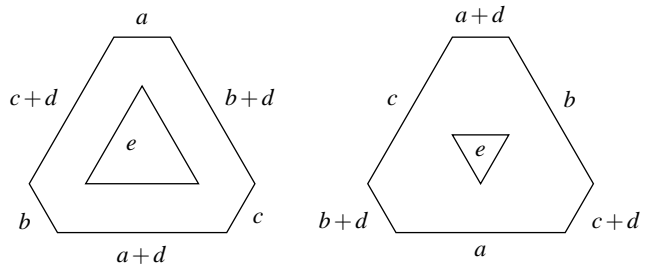
\begin{fullfigure}{f:hextriang}{Hexagonal regions without central triangles.}
\psset{xunit=.75cm,yunit=.650cm}
\pspicture(-.5,-1)(11,5)
\pspolygon(0,1)(.5,0)(3.5,0)(4,1)(2.5,4)(1.5,4)
\pspolygon(1,1)(3,1)(2,3)
\uput[30](3.25,2.5){$b+d$}\uput[90](2,4){$a$}
\uput[150](.75,2.5){$c+d$}\uput[210](.25,.5){$b$}
\uput[270](2,0){$a+d$}\uput[330](3.75,.5){$c$}
\uput[330](1.5,2){$e$}
\pspolygon(6,1)(6.5,0)(9.5,0)(10,1)(8.5,4)(7.5,4)
\pspolygon(8,1)(8.5,2)(7.5,2)
\uput[30](9.25,2.5){$b$}\uput[90](8,4){$a+d$}
\uput[150](6.75,2.5){$c$}\uput[210](6.25,.5){$b+d$}
\uput[270](8,0){$a$}\uput[330](9.75,.5){$c+d$}
\rput(8,1.667){$e$}
\endpspicture\eatline
\end{fullfigure}

As we mentioned in \sec{s:snf}, we can consider the Kasteleyn or
Kasteleyn-Percus cokernel for impossible enumerations where there are no
matchings.  This allows us to vary the above graphs in several ways, which we
describe case by case:
\begin{description}
\item[$\bullet$] We can consider $Z_\kappa(a,b,c)$ when all three of $a$, $b$,
and $c$ are odd, as well as $Z_\braket{\kappa,\rho}(a,a,a)$ when $a$ is odd.
\item[$\bullet$] If one of $a$, $b$, and $c$ is odd, we define the graph
$Z'_\kappa(a,b,c)$ by removing the central edge $e$ but not its vertices before
quotienting by $\kappa$.  If two are odd we define $Z'_\kappa(a,b,c)$ by
removing $e$ and its vertices before quotienting by $\kappa$.
\item[$\bullet$] If $Z_G(a,b,c)$ has a polygamous vertex, we can give it the
wrong parity to make $Z'_G(a,b,c)$.  If $G = \braket{\tau}$ or
$G = \braket{\rho,\tau}$, then we also define the $q$-weighted forms
$Z'_G(a,b,c;q)$ and $\tZ'_G(a,b,c;q)$.
\item[$\bullet$] We define the graph $Z(a,b,c,d,e)$ by removing a triangle of
size $e$ from a semiregular hexagon of side lengths $a,b+d,c,a+d,b,c+d$, with
$d \ne e$, as shown in \fig{f:hextriang}.  The parameters $d$ and $e$ might
even have opposite sign, which we interpret by turning the triangle
upside-down, as also indicated in \fig{f:hextriang}.
\end{description}

For each of the graphs defined in this section, we denote a corresponding
Kasteleyn-Percus matrix by replacing $Z$ by $M$ if it is bipartite and
monogamous, and the corresponding Kasteleyn matrix by replacing $Z$ by $A$
otherwise.

\begin{conjecture} Each of the matrices $M(a,b,c;q)$, $M_\rho(a,a,a;q)$,
$A_\tau(a,b,b;q)$, $\tA_\tau(a,b,b;q)$, $A'_\tau(a,a,b;q)$,
$\tA'_\tau(a,a,b;q)$, $\tA_\braket{\rho,\tau}(a,a,a;q)$,
$\tA'_\braket{\rho,\tau}(a,a,a;q)$, and $M(\mu;\vq_a)$ admits a Smith normal
form over $\Z[q,q^{-1}]$, and the entries are $q$-round.  The Smith normal form
over $\Z$ of each of the matrices $M(a,b,c,d,e)$, $M_{\kappa\tau}(2a,b,b)$,
$M_\braket{\rho,\kappa\tau}(2a,2a,2a)$, $A_\kappa(a,b,c)$, $A'_\kappa(a,b,c)$,
$A_\braket{\kappa,\tau}(2a,b,b)$, $A'_\braket{\kappa,\tau}(2a,b,b)$,
$A_\braket{\rho,\kappa}(a,a,a)$, $A_\braket{\rho,\kappa,\tau}(2a,2a,2a)$, and
$A'_\braket{\rho,\kappa,\tau}(2a,2a,2a)$ has round entries.
\label{c:round}
\end{conjecture}

We could have stated conjecture~\ref{c:round} in greater generality by
combining more of the variations above.  For example the graph $Z(a,a,a,d,e)$
has symmetries and we could consider the matrix of a suitably modified quotient
graph.  Indeed we do not know how to state Conjecture~\ref{c:round} in full
generality, since there are yet other variations of counting lozenge tilings in
a hexagon with round enumerations \cite{CK:centered}.  Presumably many or all
of these variations also have impossible counterparts, and some may have
$q$-enumerations as well.  Also the spirit of the conjecture is to find the
Smith normal forms or the cokernels explicitly.

\begin{conjecture} The Smith normal forms over $\Z[q,q^{-1}]$ of each of the
matrices $M(a,b,c;q)$, $M_\rho(a,a,a;q)$, $A_\tau(a,b,b;q)$,
$\tA_\tau(a,b,b;q)$, and $\tA_\braket{\rho,\tau}(a,a,a;q)$ are square free.
\label{c:sqfree}
\end{conjecture}

Conjecture~\ref{c:sqfree} may also not be fully general, although we note that
the Smith normal form of $M(\lambda;\vq_a)$ is not always square free. Note
that Conjectures~\ref{c:sqfree} and \ref{c:round} would together solve problem
5 in Propp's problem list \cite{Propp:progress}, which asks for the Kasteleyn
cokernel of the graph $Z(a,b,c)$.  The reason is that the Smith normal form of
any non-singular matrix $M$ is uniquely determined by $\det M$, provided that
the normal form exists and is square free. For example the $q$-enumeration of
plane partitions tells us that
$$\det M(2,2,2;q) = (2)_q^2(5)_q,$$
where
$$(n)_q = \frac{q^n-1}{q-1}.$$
So Conjectures \ref{c:round} and \ref{c:sqfree} assert that $\Sm(M(2,2,2;q))$
exists and its nontrivial entries are $(2)_q$ and $(2)_q(5)_q$.  If we
specialize at $q=1$, we obtain the correct prediction that
$$\coker M(2,2,2) \cong \Z/2\oplus \Z/10.$$

Stembridge~\cite{Stembridge:q-1} noticed a family of relations, called the
$q=-1$ phenomenon, between $-1$-enumeration of some symmetry classes of plane
partitions and ordinary enumeration of other symmetry classes. We conjecture an
extension of the phenomenon that is easier to state in terms of cokernels in
some cases and Smith normal forms in others.

\begin{conjecture} If $G$ is $\braket{1}$ or $\braket{\rho}$ and if $G' =
\braket{G,\kappa}$, then
$$\Sm(A_{G'}(a,b,c)) = \Sm(M_G(a,b,c)_{-1})^2.$$
If $G$ is $\braket{\tau}$ or $\braket{\rho,\tau}$, and if $G'$ is,
respectively, $\braket{\kappa\tau}$ or $\braket{\rho,\kappa\tau}$, then
$$\coker A_G(a,b,c)_{-1} \cong {\coker M_{G'}(a,b,c)}^{\oplus 2}.$$
If $G$ is $\braket{\tau}$ or $\braket{\rho,\tau}$ and if $G' =
\braket{G,\kappa}$, then
$$\coker A_{G'}(a,b,c) \cong \coker A'_G(a,b,c)_{-1}.$$
\label{c:q-1}
\end{conjecture}

\subsection{Jacobi-Trudi matrices}
\label{s:jt}

One conclusion of our construction for skew tableaux is a novel
determinant formula for skew Schur functions:
$$s_{\lambda/\mu}(\vx) = \det M(\lambda/\mu;\vx).$$
We can compare this to two other determinant formulas for skew
Schur functions \cite[\S7.16]{Stanley:enumerative2}. Let
$$J(\lambda/\mu;\vx)_{i,j} = h_{\lambda_i-\mu_j-i+j}(\vx),$$
where $h_n(\vx)$ is the $n$th complete symmetric function of $\vx$,
and let
$$D(\lambda/\mu;\vx)_{i,j} = e_{\lambda'_i-\mu'_j-i+j}(\vx),$$
where $e_n(\vx)$ is the $n$th elementary symmetric function of $\vx$
and $\lambda'$ is the partition conjugate to $\lambda$.
Then the Jacobi-Trudi identity states that
$$s_{\lambda/\mu}(\vx) = \det J(\lambda/\mu;\vx)$$
while the dual Jacobi-Trudi identity states that
$$s_{\lambda/\mu}(\vx) = \det D(\lambda/\mu;\vx).$$

\begin{theorem} The Jacobi-Trudi matrices $J(\lambda/\mu;\vx)$ and
$D(\lambda/\mu;\vx)$ are stably equivalent over the ring $\Z[\vx]$ to the
matrix $M(\lambda/\mu;\vx)$.
\label{th:jt}
\end{theorem}

\begin{fullfigure}{f:lmsquare}{First Gessel-Viennot model of lozenge tilings.
    (Not all orientations are shown.)}
\psset{xunit=1cm,yunit=1cm}
\pspicture(0,-.5)(3.5,4)
\multips(0,0)(1,0){4}{\psline(0,0)(0,3)}
\multips(0,0)(0,1){4}{\psline(0,0)(3,0)}
\pcline(0,3)(1,3)\mto\pcline(1,3)(1,2)\mto
\pcline(0,3)(0,2)\mto\pcline(0,2)(1,2)\mto
\pscircle[fillstyle=solid](0,3){.2}\pscircle[fillstyle=solid](2,3){.2}
\pscircle[fillstyle=solid](2,0){.2}\pscircle[fillstyle=solid](3,0){.2}
\multips(0,0)(1,0){4}{\qdisk(0,0){.1}}\multips(0,1)(1,0){4}{\qdisk(0,0){.1}}
\multips(0,2)(1,0){4}{\qdisk(0,0){.1}}\multips(0,3)(1,0){4}{\qdisk(0,0){.1}}
\rput[b](.5,3.1){$x_4$}\rput[b](1.5,3.1){$x_4$}\rput[b](2.5,3.1){$x_4$}
\rput[b](.5,2.1){$x_3$}\rput[b](1.5,2.1){$x_3$}\rput[b](2.5,2.1){$x_3$}
\rput[b](.5,1.1){$x_2$}\rput[b](1.5,1.1){$x_2$}\rput[b](2.5,1.1){$x_2$}
\rput[b](.5, .1){$x_1$}\rput[b](1.5, .1){$x_1$}\rput[b](2.5, .1){$x_1$}
\rput[t](1.5,-.5){$\lambda$}\rput[t](1.5,4){$\mu$}\rput[l](3.3,1.5){$a-1$}
\endpspicture
\end{fullfigure}

\begin{proof} The theorem is a special case of Corollary \ref{c:gv}.  Let
$X(\lambda/\mu;\vx)$ be the graph exemplified in \fig{f:lmsquare}:  a square
grid of height $a-1$ and width $\lambda_1+b$.  On the top there is a left
endpoint at each position $\mu_i+b+1-i$ and on the bottom a right endpoint at
each position $\lambda_i+b+1-i$.  Each horizontal edge points to the right and
on the $i$th row has weight $x_i$.  Each vertical edge points down.  By a
standard argument the Gessel-Viennot matrix of $X(\lambda/\mu;\vx)$ is
$J(\lambda/\mu;\vx)$. At the same time, one can check that the transit-free
resolution of $X(\lambda/\mu;\vx)$ is $Z(\lambda/\mu;\vx)$.  Thus Corollary
\ref{c:gv} says that $J(\lambda/\mu;\vx)$ is equivalent to
$M(\lambda/\mu;\vx)$.

\begin{fullfigure}{f:lmdiamond}{Second Gessel-Viennot model of lozenge tilings.
    (Not all orientations are shown.)}
\pspicture(0,-.5)(10,4.8)
\multips(0,4)(2,0){4}{\psline(0,0)(4,-4)}
\psline(1,3)(2,4)\psline(2,2)(4,4)\psline(3,1)(6,4)
\psline(4,0)(7,3)\psline(6,0)(8,2)\psline(8,0)(9,1)
\pcline(2,4)(1,3)\mto\pcline(0,4)(1,3)\mto
\pscircle[fillstyle=solid](2,4){.2}\pscircle[fillstyle=solid](6,4){.2}
\pscircle[fillstyle=solid](4,0){.2}\pscircle[fillstyle=solid](6,0){.2}
\multips(0,4)(1,-1){5}{\qdisk(0,0){.1}}
\multips(2,4)(1,-1){5}{\qdisk(0,0){.1}}
\multips(4,4)(1,-1){5}{\qdisk(0,0){.1}}
\multips(6,4)(1,-1){5}{\qdisk(0,0){.1}}
\rput[br](1.4,3.6){$x_4$}\rput[br](3.4,3.6){$x_4$}\rput[br](5.4,3.6){$x_4$}
\rput[br](2.4,2.6){$x_3$}\rput[br](4.4,2.6){$x_3$}\rput[br](6.4,2.6){$x_3$}
\rput[br](3.4,1.6){$x_2$}\rput[br](5.4,1.6){$x_2$}\rput[br](7.4,1.6){$x_2$}
\rput[br](4.4, .6){$x_1$}\rput[br](6.4, .6){$x_1$}\rput[br](8.4, .6){$x_1$}
\rput[bl](8.3,2.3){$a$}\rput[b](3,4.5){$\mu'$}\rput[t](7,-.5){$\lambda'$}
\endpspicture
\end{fullfigure}

Let $Y(\lambda/\mu;\vx)$ be the graph exemplified in \fig{f:lmdiamond}. It is a
diamond grid of height $a$ and width $\lambda_1+b$.  All edges point down, and
those in the $i$th row that point southwest have weight $x_i$.  On the top row
there is a left endpoint at each position $b+i-\mu'_i$ and on the bottom a
right endpoint at each position $b+i-\lambda'_i$.  Note that the left endpoints
are at all of the positions \emph{not} of the form $\mu_i+b+1-i$ and likewise
the right endpoints are at the positions \emph{not} of the form One can check
that the Gessel-Viennot matrix of $X(\lambda/\mu;\vx)$ is
$D(\lambda/\mu;\vx)$.  At the same time, noting the gaps between the marked
endpoints,  one can check that the transit-free resolution of
$Y(\lambda/\mu;\vx)$ is $Z(\lambda/\mu;\vx)$. Thus Corollary \ref{c:gv} says
that $J(\lambda/\mu;\vx)$ is equivalent to $D(\lambda/\mu;\vx)$.
\end{proof}

\begin{question} Are $J(\lambda/\mu;\vx)$ and $D(\lambda/\mu;\vx)$
stably equivalent over the ring of symmetric functions?
\end{question}

\thm{th:jt} implies that there are several ways that equivalent forms of the
matrix $M(\lambda/\mu;\vx)$ arise in several common guises.  If we set
$\vx=\vq_a$ and $\mu=0$, then  Conjecture~\ref{c:round} asserts that 
$M(\lambda;\vq_a)$ admits a Smith normal form  a Smith normal form over the
ring $\Z[q,q^{-1}]$. This suggests that $\coker M(\lambda;\vq_a)$ is an
important extra structure that one can associate to the representation
$V(\lambda)$ mentioned previously.  Moreover the relationship between
$M(\lambda)$ and $V(\lambda)$ is an interesting special case of
Conjecture~\ref{c:action}.

\subsection{Domino tilings}
\label{s:domino}

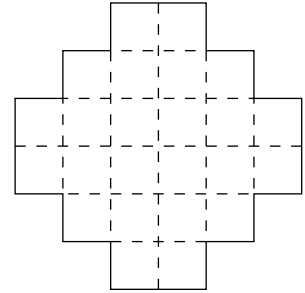
\begin{fullfigure}{f:aztec}{An Aztec diamond of order 3}
\pspicture(-3,-3)(3,3)
\pspolygon(-1,3)(1,3)(1,2)(2,2)(2,1)(3,1)(3,-1)(2,-1)(2,-2)(1,-2)(1,-3)
    (-1,-3)(-1,-2)(-2,-2)(-2,-1)(-3,-1)(-3,1)(-2,1)(-2,2)(-1,2)
\pspolygon[linestyle=dashed](-2,-1)(-2,1)(2,1)(2,-1)
\pspolygon[linestyle=dashed](-1,-2)(-1,2)(1,2)(1,-2)
\psline[linestyle=dashed](0,-3)(0,3)
\psline[linestyle=dashed](-3,0)(3,0)
\endpspicture
\end{fullfigure}

Domino tilings of an Aztec diamond are a well-known analogue of lozenge tilings
of a hexagon.  Recall that an Aztec diamond of order $n$ is the polyomino
consisting of those unit squares lying entirely inside the region $|x|+|y|
\le n+1$, as shown in \fig{f:aztec}. A domino tiling of an Aztec diamond
corresponds to a matching of the graph dual to the tiling by squares; an
example of such a graph is the one on the left in \fig{f:polytied}.  Denote
this graph by $Z_A(n)$ and let $M_A(n)$ be a Kasteleyn-Percus matrix for it.

\begin{theorem} If $M_A(n)$ is a Kasteleyn-Percus matrix for
domino tilings of an Aztec diamond of order $n$, then
$$\coker M_A(n) \cong \Z/2 \oplus \Z/4 \oplus \dots \oplus \Z/2^n.$$
\eatline
\label{th:acoker}
\end{theorem}

\thm{th:acoker} extends the result that the number of
domino tilings of an Aztec diamond of order $n$ is $2^{n(n+1)/2}$
\cite{MRR:asm,Kuperberg:eklp1,Kuperberg:eklp2}.

\begin{proof}

We give two arguments, one using Kasteleyn-Percus matrices and the other using
Gessel-Viennot matrices.  First, define a binomial coefficient matrix $B(n)$
by
$$B(n)_{0\leq i,j < n} = \binom{i}{j}.$$
(We will assume that rows and columns of other matrices are numbered from 0 as
well.) Define an $n\times n+1$ matrix $L(n)$ by putting the $n\times n$ identity
matrix to the left of a null column, and define $R(n)$ by putting the $n\times
n$ identity matrix to the right of a null column.  For example here are $B(4)$,
$L(3)$, and $R(3)$:
$$\begin{pmatrix}1 & 0 & 0 & 0\\1 & 1 & 0 & 0\\1 & 2 & 1 & 0\\1 & 3 & 3 & 1
\end{pmatrix} \quad
\begin{pmatrix}1 & 0 & 0 & 0\\0 & 1 & 0 & 0\\0 & 0 & 1 & 0 \end{pmatrix} \quad
\begin{pmatrix}0 & 1 & 0 & 0\\0 & 0 & 1 & 0\\0 & 0 & 0 & 1 \end{pmatrix}. $$

One can check that the matrix
\begin{align*}
M_A(n) =&\ R(n) \tensor R(n)^T + L(n) \tensor R(n)^T \\
&\ + R(n) \tensor L(n)^T - L(n) \tensor L(n)^T
\end{align*}
is a Kasteleyn-Percus matrix for the graph $Z_A(n)$.  At the same time,
$$B(n)L(n)B(n+1)^{-1} = L(n)$$
and
$$B(n)R(n)B(n+1)^{-1} = R(n)-L(n).$$
It follows that
\begin{align*}
M'(n) &\stackrel{\mathrm{def}}{=}
    (B(n) \tensor B(n+1)^{-T})M(n)(B(n+1)^{-1} \tensor B(n)^T) \\
    &= R(n) \tensor R(n)^T - 2L(n)\tensor L(n)^T.
\end{align*}
Since $B(n)$ is a triangular matrix, $M'(n)$ is equivalent to $M(n)$.  On
the other hand, one can check that, modulo permuting rows and columns,
\begin{align*}
M'(n) = &\ X(1)\oplus X(2)\oplus \dots \oplus X(n) \\
    &\ \oplus Y(1) \oplus Y(2) \oplus \dots \oplus Y(n),
\end{align*}
where
$$X(k)=\begin{pmatrix}1&-2&0&&0\\0&1&-2&\cdots&0\\0&0&1&&0\\&
\vdots&&\ddots&\vdots\\0&0&0&\cdots&1\end{pmatrix}$$
and
$$Y(k)=\begin{pmatrix}-2&1&0&&0\\0&-2&1&\cdots&0\\0&0&-2&&0\\&
\vdots&&\ddots&\vdots\\0&0&0&\cdots&-2\end{pmatrix}.$$
Since
$$\coker X(k)=0 \qquad \coker Y(k) =\Z/2^k,$$
the theorem follows.

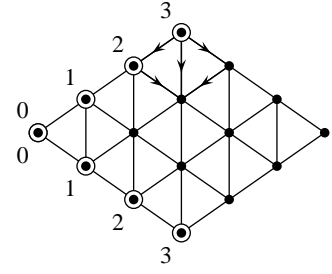
\begin{fullfigure}{f:delannoy}{Gessel-Viennot model of an Aztec diamond.
(Not all orientations are shown.)}
\pspicture(-3,-2.8)(3,2.8)
\pcline(0,2.1)(-1,1.4)\mto\pcline(0,2.1)(1,1.4)\mto
\pcline(-1,1.4)(0,0.7)\mto\pcline( 1,1.4)(0,0.7)\mto
\pcline(0,2.1)(0,0.7)\mto
\psline(0,2.1)(0,-2.1)(3,0)(0,2.1)(-3,0)(0,-2.1)
\pspolygon(-2,-0.7)(-2,0.7)(1,-1.4)(1,1.4)
\pspolygon(-1,-1.4)(-1,1.4)(2,-0.7)(2,0.7)
\multips(-3,0)(1,-0.7){4}{\pscircle[fillstyle=solid](0,0){.2}}
\multips(-3,0)(1,0.7){4}{\pscircle[fillstyle=solid](0,0){.2}\qdisk(0,0){.1}}
\multips(-2,-0.7)(1,0.7){4}{\qdisk(0,0){.1}}
\multips(-1,-1.4)(1,0.7){4}{\qdisk(0,0){.1}}
\multips(0,-2.1)(1,0.7){4}{\qdisk(0,0){.1}}
\rput[rb](-3.2,.3){$0$}\rput[rb](-2.2,1){$1$}
\rput[rb](-1.2,1.7){$2$}\rput[rb](-0.2,2.4){$3$}
\rput[rt](-3.2,-.3){$0$}\rput[rt](-2.2,-1){$1$}
\rput[rt](-1.2,-1.7){$2$}\rput[rt](-0.2,-2.4){$3$}
\endpspicture\eatline
\end{fullfigure}

In the Gessel-Viennot approach we first observe that if $G(n)$ is a graph of
the type in \fig{f:delannoy} then its transit-free resolution $G'(n)$ as
described in \thm{th:gv} is the Aztec diamond graph $Z_A(n)$.  (The left and
right endpoints of $G(n)$ are numbered in the figure, with the left endpoints
on top.  The first left and right endpoints coincide, which is degenerate but
allowed.)  Let $V(n)$ be the Gessel-Viennot matrix of $G(n)$.  The entry
$V(n)_{i,j}$ is the Delannoy number $D(i,j)$ \cite[\S6.3]{Stanley:enumerative2}
(see also Sachs and Zernitz \cite{SZ:dimer}), because $G(n)$ matches the
defining recurrence
\begin{align*}
D(i,j) = D(i,j-1) + D(i-1,j) + D(i-1,j-1) \\
D(0,i) = D(i,0) = 1.
\end{align*}
A standard formula for Delannoy numbers is
$$D(i,j)=\sum_k\binom{i}{k}\binom{j}{k}2^k.$$
In matrix form this identity can be expressed
$$V(n)=B(n)V'(n)B(n)^T,$$
where $V'(n)$ is the Smith normal form of $V(n)$ with
$$V'(n)_{k,k}=2^k.$$
Thus the Smith normal form of $V(n)$ also establishes the theorem.
\end{proof}

Finally, we could at least conjecturally extend \thm{th:acoker} with the same
variations as those we considered for lozenge tilings:  Domino tilings with
symmetry, $q$-enumerations, impossible enumerations, Aztec diamonds with teeth
missing, etc.  For example, Tokuyama \cite{Tokuyama:gelfand}
established a relation between generating functions of lozenge-type and
Aztec-type Gelfand triangles.  We conjecture that this identity can be
extended to an equivalence of Kasteleyn-Percus (or Gessel-Viennot) matrices.
We leave this and other possibilities to future work.

\section{Appendix: Smith normal form}
\label{s:smith}

\begin{theorem} If $M$ is a $k \times n$ matrix over a principal ideal domain
$R$, then there exist invertible matrices $A$ and $B$ such that
$$\Sm(M) = AMB$$
is diagonal and $\Sm(M)_{i,i}$ divides $\Sm(M)_{i+1,i+1}$.  The matrix
$\Sm(M)$ is a \emph{Smith normal form} of $M$.
\label{th:smith}
\end{theorem}

Evidently
\begin{align*}
\coker M &\cong \coker \Sm(M) \\
    &= R^{k-n} \oplus R/\Sm(M)_{1,1} \oplus \dots \oplus R/\Sm(M)_{n,n}
\end{align*}
if $k \ge n$ and
\begin{align*}
\coker M &\cong \coker \Sm(M) \\
    &= R/\Sm(M)_{1,1} \oplus \dots \oplus R/\Sm(M)_{k,k}
\end{align*}
otherwise.  (If $R = \Z$ and we allow $n$, but not $k$, to be infinite, then
\thm{th:smith} is equivalent to the classification of finitely generated
abelian groups.)  It is easy to show that this decomposition of $\coker M$ is
unique and that $\Sm(M)$ is unique up to multiplying the entries by units.

\begin{proof}
We argue by induction on $n$ or $k$. Since $R$ is a PID, any two elements $a$
and $b$ have a greatest common divisor
$$\gcd(a,b) = ax+by$$
which is unique up to a unit factor.  Moreover, $R$ is Noetherian, which means
that any chain of divisibilities $a_{i+1}|a_i$ must eventually be constant up
to unit factors.  Thus we can argue by induction with respect to divisibility
of non-zero elements of $R$.   Since we can replace $M$ by any equivalent form,
we assume that the entry $M_{1,1} = a$ is not divisible by any entry in any
form of $M$, except for those that are $a$ times a unit.  We claim that $a$
divides every entry of $M$.  Otherwise there is an entry $M_{i,j} = b$ such
that $a \nmid b$ and $b \nmid a$.  If $b$ is in the first row or column, then
after permuting rows and columns and possibly transposing $M$, $M_{1,2} = b$:
$$M = \left(\begin{array}{cc|c}
    a & b & \cdots \\ \hline \multicolumn{2}{c|}{\vdots} & \ddots
\end{array}\right).$$
Let $c = ax+by$ be a greatest common divisor of $a$ and $b$. If we
post-multiply $M$ by the matrix
$$B = \left(\begin{array}{c|c}
    \begin{matrix}x & a/c \\ y & -b/c \end{matrix} & \mbox{\Large $0$}
    \\ \hline \mbox{\Large \rule{0ex}{2ex} $0$} & \mbox{\Large $I$}
\end{array}\right),$$
the result is a form of $M$ with the entry $c$, which contradicts the choice of
$a$.  If $a$ divides every entry in the first row and column but not some entry
$b$ elsewhere, then we first pivot at the position $(1,1)$ to obtain:
$$M = \left(\begin{array}{cc|c}
    a & 0 & \ldots \\ 0 & b & \\ \hline \multicolumn{2}{c|}{\vdots} & \ddots
\end{array}\right).$$
If we pre-multiply by
$$A = \left(\begin{array}{c|c}
    \begin{matrix}1 & 1 \\ 0 & 1 \end{matrix} & \mbox{\Large $0$}
    \\ \hline \mbox{\Large \rule{0ex}{2ex} $0$} & \mbox{\Large $I$}
\end{array}\right)$$
and post-multiply by $B$ above,
we again produce the entry $c$.

Given that $M_{1,1} = a$ does divide the rest of $M$ we perform
a deleted pivot at $(1,1)$ and assume normal form for the remaining
submatrix by induction.
\end{proof}

\begin{theorem} If $A$ is an alternating $n \times n$ matrix over a principal
ideal domain $R$, then there exists an invertible matrix $B$ such that
$$\Sm_a(A) = B^TAB$$
is block-diagonal with $2 \times 2$ blocks, and such that $\Sm_a(A)_{2i,2i+1}$
divides $\Sm_a(A)_{2i+2,2i+3}$.  The matrix $\Sm_a(A)$ is the
\emph{alternating Smith normal form} of $A$.
\label{th:asmith}
\end{theorem}
\begin{proof} The argument is the same as the one for Theorem~\ref{th:smith}.
We assume that $A_{1,2} = a$ is minimal among all entries of all forms of $A$.
If it does not divide some entry in the first two rows and columns, then
after permuting rows and columns, that entry is $A_{1,3} = b$:
$$A = \left(\begin{array}{ccc|c}
    0 & a & b &  \\ -a & 0 & \cdot & \cdots \\
    -b & \cdot & 0 & \\ \hline & \vdots & & \ddots
\end{array}\right).$$
Let $c = ax+by$ be a common divisor and let
$$B = \left(\begin{array}{c|c}
    \begin{matrix}1 & 0 & 0  \\ 0 & x & a/c \\ 0 & y & -b/c \end{matrix}
    & \mbox{\Large $0$} \\ \hline \mbox{\Large \rule{0ex}{2ex} $0$}
    & \mbox{\Large $I$}
\end{array}\right).$$
Then $B^TAB$ has the entry $c$, a contradiction.  If $a$ divides every entry in
the first two rows and columns but not some other entry $b$, then we can
perform a symmetric pivot at $(1,2)$.  After permuting rows and columns $A$
then has the form:
$$A' = \left(\begin{array}{cccc|c}
      0 & a &  0 & 0 & \\
     -a & 0 &  0 & 0 & \ldots \\
      0 & 0 &  0 & b & \\
      0 & 0 & -b & 0 & \\ \hline
      \multicolumn{4}{c|}{\vdots} & \ddots
\end{array}\right).$$
Let $c = ax+by$ be a common divisor and let
$$B = \left(\begin{array}{c|c}
    \begin{matrix}
    1 & 0 & 0 & 0 \\
    0 & 1 & 0 & 0 \\
    1 & 0 & 1 & 0 \\
    0 & 0 & 0 & 1 \end{matrix}
    & \mbox{\Large $0$} \\ \hline \mbox{\Large \rule{0ex}{2ex} $0$}
    & \mbox{\Large $I$}
\end{array}\right).$$
Then $B^TAB$ has the form of the previous case.

Finally if $a$ divides every entry of $A$, we perform a deleted symmetric pivot
at $(1,2)$ and inductively assume the normal form for the remaining submatrix.
\end{proof}


\providecommand{\bysame}{\leavevmode\hbox to3em{\hrulefill}\thinspace}

\end{document}